\documentclass[12pt]{article}
\usepackage[landscape]{geometry}
\usepackage[francais]{babel}

\newtheorem{theoreme}{Th\'eor\`eme}
\newtheorem{proposition}{Proposition} 
\newtheorem{definition}{D\'efinition} 
\newtheorem{corollaire}{Corollaire}  

\newtheorem{lemme}{Lemme} 
 
\newtheorem{remarque}{Remarque}
\newtheorem{res}{R\'esultat}

 \title{Filtration  Relative, l'Id\'eal de Bernstein    et  ses pentes}
 
 \author{ \vspace{1cm} \\ Ph. Maisonobe \\
   Universit\'e de Nice Sophia Antipolis \\
   Laboratoire Jean-Alexandre Dieudonn\'e \\
  Unit\'e Mixte de Recherche du  CNRS  7351 \\ 
  Parc Valrose, F-06108 Nice Cedex 2}

\begin{document}
\maketitle

 \newpage

\section*{Introduction   }

Soit $f_i : X \rightarrow  {\bf C}$, pour $i$ entier compris entre $1$ et $p$, des fonctions analytiques d\'efinies   au voisinage d'un compact $K$
d'une  vari\'et\'e analytique complexe $X$. Notons   $F$ le produit des $f_i$ et posons si     $\phi  : X \rightarrow  {\bf C}$ d\'esigne une fonction $C^{\infty}$ \`a support compact dans $K$   :
$$ I_{\phi} (s_1, \ldots,s_p) = \int_X \mid f_1 (x) \mid^{s_1}\ldots \mid f_p (x) \mid^{s_p} \phi (x) \;  dx \wedge d\overline{x} \; .$$
Tout comme dans le cas $p=1$, en utilisant le th\'eor\`eme de r\'esolution des singularit\'es d'H. Hironaka, on peut montrer que $ I_{\phi} (s_1, \ldots,s_p)$ qui est une fonction d\'efinie \`a priori
pour ${\rm Re} \, s_i  > 0$ se prolonge en fonction m\'eromorphe avec des p\^oles situ\'es sur des hyperplans de ${\bf C}^p$ (th\'eor\`eme 1 de \cite{K-K}). 
F. Loeser \'etudie ces int\'egrales dans \cite{L} et    appelle pente de $(f_1, \ldots ,f_p)$  les directions de leurs  hyperplans polaires.   Dans
certains cas g\'eom\'etriques, il majore  cet ensemble de  pentes par un ensemble de  formes lin\'eaires 
li\'ees  \`a la g\'eom\'etrie du discriminant du morphisme $(f_1, \ldots ,f_p) : X \rightarrow  {\bf C}^p$.\\

Consid\'erons  ${\cal D}_X$ l'anneau des op\'erateurs diff\'erentiels et    ${\cal D}_X[s_1, \ldots ,s_p]= {\bf C}_X[s_1, \ldots ,s_p] \otimes_{\bf C} {\cal D}_X$.
Soit $m$ une section d'un ${\cal D}_X$-Module holonome,  notons ${\cal B}(m,x_0, f_1, \ldots ,f_p)$ l'id\'eal de   ${\bf C}[s_1, \ldots ,s_p]$   des polyn\^omes   $b$ 
  v\'erifiant  au voisinage de $x_0$:
$$    b (s_1, \ldots ,s_p) m f_1^{s_1} \ldots  f_p^{s_p} \in {\cal D}_X[s_1, \ldots ,s_p]\,  m f_1^{s_1+1} \ldots  f_p^{s_p+1} \; .$$
 Ces  polyn\^omes sont appel\'es  polyn\^omes de Bernstein de $(m, f_1, \ldots ,f_p) $ au voisinage de $x_0$. Suivant
   J. Bernstein \cite{B},  ils  permettent de construire  
 un prolongement des int\'egrales $ I_{\phi} (s_1, \ldots,s_p)$. Dans \cite{S1}, {\bf C. Sabbah montre     l'existence 
  pour tout $x_0 \in X$  d'un ensemble   fini ${\cal H}$  de formes lin\'eaires  \`a coefficients premiers entre eux dans ${\bf N}$ telles que :
    $$\prod_{H\in {\cal H} } \prod_{i\in I_{\cal H} } (H(s_1, \ldots, s_p) + \alpha _{H,i}) \in {\cal B}(m,x_0, f_1, \ldots ,f_p) \; ,$$
    o\`u $ \alpha _{H,i}$ sont des nombres complexes.  }
     Dans \cite{S2}, il  montre comment en d\'eduire   des r\'esultats analogues \` a ceux de  F. Loeser. 
   Mais,  J. Brian\c{c}on et H. Maynadier montrent dans \cite{B.May} que l'id\'eal ${\cal B}(m,x_0, f_1, \ldots ,f_p)$ n'est en g\'en\'eral pas principal.
   
   {\bf L'objet de cet article est notamment de montrer
   l'existence d'un ensemble    ${\cal H}$  minimal.}     De plus,  lorsque $m$ est une section d'un module holome r\'egulier,
 nous expliciterons  cet ensemble  g\'eom\'etriquement 
  \`a partir de la vari\'et\'e caract\'eristique du syst\`eme 
diff\'erentiel engendr\'e par $m$.   \\

  Sur  ${\cal D}_X[s_1, \ldots ,s_p]$,
 nous consid\'erons   la filtration di\`ese (resp. la filtration relative) qui \'etend la filtration de ${\cal D}_X$   en donnant \`a $s_i$  le   poids un (resp. z\'ero). 
Si $M$ est un ${\cal D}_X[s_1, \ldots ,s_p] $-Module coh\'erent muni d'une bonne filtration di\`ese (resp. relative), nous notons ${\rm gr }^{\sharp}\,  M$  (resp. ${\rm gr }^{\rm rel}\,  M$ )
son gradu\'e pour cette filtration. La racine de l'annulateur de ${\rm gr }^{\sharp}\,  M$  (resp. ${\rm gr }^{\rm rel}\,  M$ ) est ind\'ependante de la bonne filtration et  d\'efinit un  sous-espace  analytique de $T^{\ast}X \times {\bf C}^p$
appel\'e vari\'et\'e caract\'eristique di\`ese (resp. relative) de $M$ et not\'ee ${\rm car}^{\sharp} \, M$
(resp. ${\rm car}^{\rm rel} \, M)$.\\

Si $M$ est un ${\cal D}_X[s_1, \ldots ,s_p]$-Module coh\'erent, nous expliquons que les dimensions des  vari\'et\'es caract\'eristiques di\`ese et relative de $M$ 
 se d\'eterminent \` a l'aide des nombres grade des fibres de $M$.  Ces r\'esultats g\'en\'eralisent ceux  des ${\cal D}_X$-Modules coh\'erents (voir  pour ce cas  \cite{K2},   \cite{Bj2}, \cite{G-M}   ... ).
En utilisant le th\'eor\`eme d'involutivit\'e de O. Gabber, nous montrons alors que pour toute section $m$ d'un Module holonome :\\

\begin{res}
  Il existe une vari\'et\'e 
  lagrangienne conique $\Lambda$  de $T^{\ast}X$ non lisse en g\'en\'eral  telle que :
    $$  {\rm car }^{\rm rel}\,  {\cal D}_X[s_1, \ldots ,s_p] m f_1^{s_1} \ldots  f_p^{s_p}   =  \Lambda \times  {\bf C}^p \; .\\$$
  \end{res}
  
 Il en r\'esulte en particulier que la dimension du  
 ${\cal D}_X[s_1, \ldots ,s_p]$-Module 
 $$\frac{ {\cal D}_X[s_1, \ldots ,s_p] m f_1^{s_1} \ldots  f_p^{s_p} }{{\cal D}_X[s_1, \ldots ,s_p] m f_1^{s_1+1} \ldots  f_p^{s_p+1}} $$
 est    inf\'erieure  ou \'egale \`a $ {\rm dim}\, X +p -1$.  
 
\begin{res} Il existe  une famille de couples $( S_{\alpha},X_{\alpha})$ o\`u les $ S_{\alpha}$  sont des
  sous-vari\'et\'es  alg\'ebriques    de  ${\bf C}^p$ de dimensions  inf\'erieures  ou \'egales  \`a $p-1$ et les   $X_{\alpha}$ des   sous-espaces  analytiques de  $X$
  telles que~:
 $$  {\rm car}^{\rm rel} \,  \frac{ {\cal D}_X[s_1, \ldots ,s_p] m f_1^{s_1} \ldots  f_p^{s_p} }{{\cal D}_X[s_1, \ldots ,s_p] m f_1^{s_1+1} \ldots  f_p^{s_p+1}} = 
 \bigcup\,   T_{X_{\alpha} }^{\ast} X  \times  S_{\alpha}   \; .$$
   \end{res}

Par d\'efinition de la  vari\'et\'e caract\'eristique relative, {\bf la r\'eunion des $S_{\alpha}$  pour $x_0  \in  X_{\alpha}$
 n'est autre que la vari\'et\'e des z\'eros de l'id\'eal
 ${\cal B}(m,x_0, f_1, \ldots ,f_p)$}. Au passage,  pour  tout  $x_0$ dans $X$, 
notre r\'esultat  donne    une nouvelle d\'emonstration     de l'existence d'un polyn\^ome non nul dans ${\cal B}(m,x_0, f_1, \ldots ,f_p)$. 
Cette d\'emonstration n'est int\'eressante que dans le cas analytique En effet,
si les   $f_i$ sont des polyn\^omes,  
 la preuve de J. Bersntein donn\'ee dans \cite{B} se g\'en\'eralise    sans modification et reste la meilleure r\'ef\'erence. \\
 
Pour pr\'eciser la structure de l'id\'eal  ${\cal B}(m,x_0, f_1, \ldots ,f_p)$, le probl\`eme est que  les vari\'et\'es alg\'ebriques   $S_{\alpha}$ ne sont \`a priori pas toutes r\'eunions d'hypersurfaces. 
 En adaptant un r\'esultat   inspir\'e de O. Gabber  et donn\'e par Y.E.  Bj\"{o}rk  dans \cite{Bj2}  sur le   conoyau d'un endomorphisme injectif d'un Module  pur, nous montrons que si 
   $M$ est un ${\cal D}_X$-Module  holonome    engendr\'e par une section $m$, 
la
 famille ${\cal G}$  des     sous-$ {\cal D}_{X}[s_1, \ldots ,s_p]$-Modules  $L$ de  type fini de $\displaystyle M[\frac{1}{F},s_1, \ldots ,s_p] f_1^{s_1 } \ldots  f_p^{s_p } $   contenant  $ {\cal D}_{X}[s_1, \ldots ,s_p] m f_1^{s_1 } \ldots  f_p^{s_p } $
et  tels que pour  tout point $x_0 \in X$ :
$$  {\rm grade} \, \frac{L_{x_0}}{ {\cal D}_{X,x_0}[s_1, \ldots ,s_p] m f_1^{s_1 } \ldots  f_p^{s_p } } \geq {\rm dim}\, X   +  2 \; $$
 admet un plus grand \'el\'ement  not\'e   $\tilde{L}$. Ce Module $\tilde{L}$ est   un 
 ${\cal D}_{X}[s_1, \ldots ,s_p]$-Module coh\'erent   
 v\'erifiant~:
\begin{enumerate}
\item ${\cal D}_{X}[s_1, \ldots ,s_p] m  f_1^{s_1} \ldots  f_p^{s_p} \subset \tilde{L}$,
\item $\tau ( \tilde{L} ) \subset \tilde{L} $,
\item $ \tilde{L} / \tau ( \tilde{L} )$ est un ${\cal D}_{X}[s_1, \ldots ,s_p]$-Module coh\'erent dont les fibres non nulles sont des modules purs de grade
${\rm dim}\, X +1$, c'est \` a dire    des modules tels que tous leurs sous-modules non r\'eduit  \`a z\'ero aient m\^eme dimension.\\
\end{enumerate}
Nous en d\'eduisons la pr\'ecision suivante : 

\begin{res} 
 $$  \displaystyle {\rm car}^{\rm rel} \,  \frac{ {\cal D}_X[s_1, \ldots ,s_p] m f_1^{s_1} \ldots  f_p^{s_p} }{{\cal D}_X[s_1, \ldots ,s_p] m f_1^{s_1+1} \ldots  f_p^{s_p+1}} = 
 \bigcup\,   T_{X_{\alpha} }^{\ast} X  \times  S_{\alpha}  \quad {\rm avec}  $$
  
 \begin{itemize}  
 \item Chaque vari\'et\'e alg\'ebrique  $S_{\alpha}$ est de dimension $p-1$.  
 \item Les composantes irr\'eductibles de 
 dimension $p-1$ de chaque $S_{\alpha}$  sont des hyperplans affines  $H_{\alpha,\beta}$  dont les directions sont  des noyaux de  formes lin\'eaires \`a coefficients entiers positifs et  premiers entre eux dans ${\bf N}$.
 \item   Les   composantes  irr\'eductibles   des  $S_{\alpha}$  de    dimension strictement inf\'erieure \`a $p-1$ sont   contenues  dans des  hyperplan affine du type $\tau ^k ( H_{\alpha,\beta}  )$ o\`u $k \in {\bf Z}$ et 
  $\tau$ la translation $(s_1, \ldots, s_p) \mapsto (s_1+1, \ldots, s_p+1)$. \\
  \end{itemize} 
  Pour $x_0 \in X $, {\bf notons ${\cal H}(x_0,m)$   l'ensemble des directions des hyperplans $H_{\alpha,\beta}$} pour $\alpha$ tel  que $x_0$ soit dans $X_{\alpha}$. Nous 
 appelons ces directions les  pentes de $(m,f_1, \dots,f_p)$ au voisinage de $x_0$. 
   \end{res}

Le fait que les composantes irr\'eductibles de dimension $p-1$ des $S_{\alpha}$ soient  des hyperplans affines dont les  directions sont les  noyaux de formes lin\'eaires \`a coefficients premiers entre eux dans ${\bf N}$ 
se d\'eduit en fait du r\'esultat de C. Sabbah sur l'existence d'un polyn\^ome de Bernstein qui soit un produit de formes lin\'eaires. \\

Nous obtenons alors :
 
\begin{res}  Condid\'erons un produit de formes 
  lin\'eaires affines appartient \`a $  {\cal B}(m,x_0, f_1, \ldots ,f_p)$. Alors, tout hyperplan vectoriel  de  ${\cal H}(x_0,m)$ est direction de l'un de ses facteurs.
De plus, il existe dans $  {\cal B}(m,x_0, f_1, \ldots ,f_p)$ un produit de formes lin\'eaires affines dont les   directions sont  exactement l'ensemble ${\cal H}(x_0,m)$ des 
    pentes de  $(m,f_1, \dots,f_p)$ au voisinage de $x_0$. Enfin, 
 la vari\'et\'e  des z\'eros de ${\bf in }\,  {\cal B}(m,x_0, f_1, \ldots ,f_p)$  l'id\'eal   engendr\'e par les  parties homog\'enes de plus haut degr\'e  
     des \'el\'ements de $ {\cal B}(m,x_0, f_1, \ldots ,f_p)$ est la r\'eunion  des   pentes de  $(m,f_1, \dots,f_p)$ au voisinage de $x_0$.
      La  racine de ${\bf in }\,  {\cal B}(m,x_0, f_1, \ldots ,f_p)$  est en particulier un    id\'eal principal.\\
\end{res}
     
     Consid\'erons l'application :
     $$ {\rm exp}^{ 2i \pi . } : {\bf C}^p \longrightarrow  ({\bf C}^{\ast})^p \quad , \quad (s_1, \ldots ,s_p) \longmapsto (e^{2i\pi s_1}, \ldots ,e^{2i\pi s_p})\; .\\$$
     
 \begin{res}     
     L'image par l'application ${\rm exp}^{ 2i \pi . }$ de la  vari\'et\'es des z\'eros de $  {\cal B}(m,x_0, f_1, \ldots ,f_p)$ est une r\'eunion
        de sous-ensembles de $ ({\bf C}^{\ast})^p$  o\`u  chaque sous-ensemble est d\'efini  par une \'equation du  type :
  $$ (\sigma _1)^{a_1} \cdots (\sigma _p)^{a_p} = \alpha $$ 
  o\` u $(a_1, \ldots ,a_p)$ est une famille d'\'el\'ements de ${\bf N}$ premier entre eux et $\alpha $ un nombre complexe. L'ensemble des $(a_1, \ldots ,a_p)$ 
   est l'ensemble des coefficients des d'\'equations  des   pentes de  $(m,f_1, \dots,f_p)$ au voisinage de $x_0$.  
   \end{res}
   
   Ce r\'esultat r\'epond \`a une  question de N. Budur \cite{Bu} pos\'ee pour le cas particulier
  $M= {\cal O}_X$.  \\

  Nous allons maintenant pr\'eciser  ${\cal H}(x_0,m)$ lorsque $m$ est une section d'un Module holonome r\'egulier.  Rappelons pour cela quelques notations.\\
  
  Soit  $\Lambda$ une vari\'et\'e lagrangienne conique de $T^{\ast}X$.
     Nous
d\'esignons par    $W^{\sharp}_{f_1 , \ldots f_p , \Lambda}$ {\bf l'adh\'erence dans $T^{\ast}X \times {\bf C}^p$ de
$$\{ (x, \xi + \sum_{i=1}^p s_i \frac{ df_i(x)}{f_i(x)}, s_1, \ldots ,s_p ) \; ; \; s_i \in {\bf C} \; , \;  (x, \xi) \in  \Lambda 
\; {\rm et} \; F(x) \neq 0 \}\; .$$}
Dans \cite{B-M-M1}, avec J. Brian\c{c}on et M. Merle, nous avions montr\'e que si $\Lambda$ n'est pas contenu dans $F^{-1}(0)$, 
les composantes irr\'eductibles de $ W^{\sharp}_{f_1, \ldots ,f_p, \Lambda} \cap F^{-1}(0)$ sont toutes de dimension $ {\rm dim}\, X +p-1$. Leurs projections sur ${\bf C}^p$  
 sont des  hyperplans vectoriels   dont les \'equations sont des formes lin\'eaires \`a coefficients entiers positifs ou nuls. 
{\bf Nous appelons pentes de $(\Lambda, f_1, \ldots ,f_p )$ au voisinage de $x_0$ les hyperplans vectoriels obtenus par projection  sur ${\bf C}^p$   des composantes irr\'eductibles de  
$ W^{\sharp}_{f_1, \ldots ,f_p, \Lambda} \cap F^{-1}(0)$ qui rencontre la fibre de $x_0$. 
Nous notons ${\cal H}(\Lambda,x_0,f_1, \ldots ,f_p )$  l'ensemble de ces pentes.\\}

  A l'aide d'un r\'esultat  de C. Sabbah sur les vari\'et\'es caract\'eristiques  d'un module relatif engendrant un module holonome 
r\'egulier (th\'eor\`eme 3.2, \cite{S2}), nous avions   \'etabli  les r\'esultats suivants avec J. Brian\c{c}on et M. Merle dans \cite{B-M-M3}.
 Soit $M$ un 
 ${\cal D}_X$-Module holonome r\'egulier $M$ de vari\'et\'e caract\'eristique $\Lambda$
  et $m$ une section   engendrant $M$,   alors~:  
  \begin{enumerate}
\item  $\displaystyle {\rm car}^{\sharp} \, ( {\cal D}_X[s_1, \ldots ,s_p]  m   f_1^{s_1} \ldots  f_p^{s_p} )  =  W^{\sharp}_{f_1, \ldots ,f_p, \Lambda}  \; ,$ 
\item $\displaystyle {\rm car}^{\sharp}\,  \left( \frac{{\cal D}_X[s_1, \ldots ,s_p]  m   f_1^{s_1} \ldots  f_p^{s_p}}{{\cal D}_X[s_1, \ldots ,s_p]  m   f_1^{s_1+1} \ldots  f_p^{s_p+1}} \right) 
=  W^{\sharp}_{f_1, \ldots ,f_p, \Lambda}  \bigcap F^{-1}(0)\; .\\$
\item $ W^{\sharp}_{f_1, \ldots ,f_p, \Lambda} (0)  = W^{\sharp}_{f_1, \ldots ,f_p, \Lambda}  \bigcap (s_1 = \ldots s_p =0)$ est une vari\'et\'e lagrangienne conique de $T^{\ast}X$.
 \end{enumerate}
 
 Un premier r\'esultat dans ce sens \'etait donn\'e par M. Kashiwara et T. Kawai dans \cite{K-K}.\\

 Nous montrons :\\
 
  \begin{res} 
Si $m$ est une section d'un Module holonome r\'egulier engendrant un  ${\cal D}_X$-Module  de vari\'et\'e caract\'eristique $\Lambda$ :
  \begin{itemize}
\item
   Les pentes   de $(m,f_1, \dots,f_p)$ au voisinage de $x_0$ sont \'egales aux pentes de  $(\Lambda, f_1, \ldots ,f_p )$ au voisinage de $x_0$ :
  $ {\cal H}(x_0,m) =  {\cal H}(\Lambda,x_0,f_1, \ldots ,f_p ) \; .  $  
\item   $  {\rm car }^{\rm rel}\,  {\cal D}_X[s_1, \ldots ,s_p] m f_1^{s_1} \ldots  f_p^{s_p}   =  W^{\sharp}_{f_1, \ldots ,f_p, \Lambda} (0) \times  {\bf C}^p \; .\\$
 \end{itemize}
  \end{res}

 Nous avions obtenu   avec J. Brian\c{c}on et M. Merle  dans \cite{B-M-M2} des r\'esultats analogues   pour $p=2$  et  
 pour $p$ quelconque dans  le cas   cas  o\`u les pentes de $(\Lambda, f_1, \ldots ,f_p )$ sont contenues dans les hyperplans de corodonn\'ees de ${\bf C}^p$. 
 Dans \cite{Mai},  est  d\'evelopp\'e une th\'eorie des cycles \'evanescents den ces  morphismes que nous appelons sans pente.

       \newpage

 \tableofcontents
\newpage

\section{ Filtration di\`ese  et filtration relative d'un ${\cal D}_X [s_1, \ldots ,s_p]$-Module}
\label{s1}
Soit $X$ une vari\'et\'e analytique  complexe. Nous d\'esignons  par  ${\cal O}_X$ le faisceau des fonctions holomorphes sur $X$
et par ${\cal D}_X$ celui des op\'erateurs diff\'erentiels    muni de sa filtration naturelle 
$( {\cal D}_X(k) )_{k \in {\bf N}}$ d\'efinie par l'ordre des d\'erivations.\\

Localement, nous identifierons $X$ \`a ${\bf C}^n$ au moyen d'un syst\`eme  $(x_1, \ldots ,x_n)$ de coordonn\'ees  locales.
Un op\'erateur $P \in  {\cal D}_X(k)  $  s'\'ecrit alors~:
$$ P = \sum_{ \mid \beta \mid \leq k}  c_{\beta}(x) \partial ^{\beta} 
\quad {\rm  o\grave{u} } \quad \beta = (\beta _1, \ldots , \beta _n  )\in {\bf N}^n 
\quad {\rm  o\grave{u} } \quad c_{\beta} \in {\cal O}_X
 \;   {\rm  et}   \; \partial ^{\beta} = \frac{ \partial ^{\beta _1} }{ \partial x_1 }  \cdots \frac{ \partial ^{\beta _n} }{\partial x_n}  \; .$$

D\'esignons  
$\pi : T^{\ast}X   \rightarrow X$ le fibr\'e cotangent \` a $X$  . Le gradu\'e ${\rm gr } \, {\cal D}_X$ s'identifie au sous-faisceau 
 de $\pi _{\ast} ( {\cal O}_{T^{\ast}X} )$ des fonctions analytiques sur $T^{\ast} X$ polynomiales par rapport aux fibres de $\pi$. ll est isomorphe 
localement au faisceau d'anneaux commutatifs  gradu\'es $  {\cal O}_ X [\xi _1  \ldots \xi _n ] $.
Si $P \in {\cal D}_X(k)  - {\cal D}_X(k-1) $, sa classe modulo ${\cal D}_X(k-1)$ d\'efinit une section de ${\rm gr }\;  {\cal D}_X$
appel\'e symbole principal de $P$ et not\'e $\sigma (P)$.\\

Si $M$ est un ${\cal D}_X$-Module coh\'erent, nous notons ${\rm car}_{{\cal D}_X}\,   M $ sa vari\'et\'e caract\'eristique.\\

Soit $p$ un entier sup\'erieur ou \'egal \` a $1$.   Dans cette section,  nous allons \'etudier le faisceau d'anneaux  ${\cal D}_X[s_1, \ldots ,s_p]= {\bf C}_X[s_1, \ldots ,s_p] \otimes_{\bf C} {\cal D}_X$.
Nous notons toujours  $\pi : T^{\ast}X \times {\bf C}^p  \rightarrow X$ la  projection naturelle.

\subsection{Filtration di\`ese de ${\cal D}_X[s_1, \ldots ,s_p] $}
 
 La filtration $( {\cal D}^{\sharp}_X[s_1, \ldots ,s_p](k) )_{k \in {\bf N}}$ 
de ${\cal D}_X[s_1, \ldots ,s_p] $ est 
d\'efinie comme suit : un op\'erateur $P$ appartient \`a  $ {\cal D}^{\sharp}_X[s_1, \ldots ,s_p](k) $ s'il s'\'ecrit localement~:
$$  P = \sum_{ \mid \alpha \mid + \mid \beta \mid \leq k}  a_{{\alpha},\beta}(x) s^{\alpha}\partial ^{\beta} \quad {\rm  o\grave{u} } 
\quad  a_{{\alpha},\beta}\in {\cal O}_X \; , \; 
s^\alpha = s_1^{\alpha _1}  \cdots s_p^{\alpha _p}
   \;,$$
Nous l'appelons cette filtration  la \underline{filtration di\`ese} de ${\cal D}_X[s_1, \ldots ,s_p] $.\\
 
Nous d\'esignons par ${\rm gr }^{\sharp} \, {\cal D}_X[s_1, \ldots ,s_p]$ le gradu\'e de cette filtration. 
Ce gradu\'e  s'identifie au sous-faisceau 
de $\pi _{\ast} ( {\cal O}_{T^{\ast}X \times {\bf C}^p} )$ des fonctions analytiques sur $T^{\ast} X \times {\bf C}^p$ polynomiales par rapport aux fibres de $\pi$. Il 
est isomorphe localement   au faisceau d'anneaux commutatifs  gradu\'es $  {\cal O}_ X [s_1, \ldots ,s_p, \xi _1,  \ldots ,\xi _n ] $. Il en r\'esulte
que ${\cal D}_X[s_1, \ldots ,s_p] $ est un faisceau coh\'erent d'anneaux 
(voir par exemple  \cite{Bj2} Appendice III, theorem 2.7).\\

  Un op\'erateur $P \in{\cal D}^{\sharp}_X [s_1, \ldots ,s_p](k)  - {\cal D}^{\sharp}_X[s_1, \ldots ,s_p](k-1)$   modulo $ {\cal D}^{\sharp}_X[s_1, \ldots ,s_p](k-1)$  d\'efinit une section de ${\rm gr }^{\sharp} {\cal D}_X[s_1, \ldots ,s_p]$
appel\'ee symbole di\`ese   de $P$ et not\'e $\sigma ^{\sharp}(P)$. Localement~: \\
$$ {\rm  si  }   \quad P  =  \sum_{ \mid \alpha \mid + \mid \beta \mid \leq k}  a_{{\alpha},\beta}(x) s^{\alpha}\partial ^{\beta} 
\quad  , \quad  
 \sigma ^{\sharp}(P) (x,\xi,s) = \sum_{ \mid \alpha \mid + \mid \beta \mid = k}  a_{{\alpha},\beta}(x) s^{\alpha}\xi ^{\beta} \; .$$

 Soit $M$  un  ${\cal D}_X[s_1, \ldots ,s_p] $-Module coh\'erent et $(M_k)_{k \in {\bf Z} }$ une filtration de $M$ pour la filtration di\`ese. Cette filtration est dite bonne si
 localement, il existe des sections $m_1, \ldots ,m_l$ de $M$  et des entiers relatifs $k_1, \ldots , k_l$ tels que $M_k = \sum_{i=1}^l {\cal D}^{\sharp}_X[s_1, \ldots ,s_p](k-k_i)   m_i$. 
 Nous notons ${\rm gr }^{\sharp}\,  M$ le gradu\'e de $M$ pour cette filtration et ${\rm Ann}\, {\rm gr }^{\sharp}\,  M$ son annutateur comme 
${\rm gr }^{\sharp} \, {\cal D}_X[s_1, \ldots ,s_p]$-Module. La racine de ${\rm Ann}\, {\rm gr }^{\sharp}\, M$ est un id\'eal coh\'erent
de ${\rm gr }^{\sharp} \, {\cal D}_X[s_1, \ldots ,s_p]$ ind\'ependante des bonnes filtrations di\`eses de $M$. Nous notons $J^{\sharp} (M)$ cet id\'eal.
Nous appelons  \underline{vari\'et\'e caract\'eristique di\`ese} de $M$  
le  sous-espace analytique de $T^{\ast}X \times {\bf C}^p  $ d\'efini par $J^{\sharp} (M)$. Nous le    notons  ${\rm car}^{\sharp} \, M$ .\\

Soit $x_0  \in X$. Au voisinage de $x_0$, nous identifions $X$ \`a ${\bf C}^n$ et $x_0$ \`a l'origine au moyen d'un syst\`eme de coordonn\'ees locales.
Pour tout $(a,b) \in {\bf C}^n \times {\bf C}^p$, notons ${\cal M}_{x_0,a,b}$ l'id\'eal maximal de ${\rm gr }^{\sharp} \, {\cal D}_{X  }[s_1, \ldots ,s_p]_{x_0}$
engendr\'e par $  x_1, \ldots ,x_n, \xi _1 - a_1, \ldots , \xi _n - a_n, s_1 - b_1,  \ldots , s_p - b_p  $. Le localis\'e 
$(({\rm gr }^{\sharp} \, M)_{x_0})_{{\cal M}_{x_0,a,b}} $ est un module de type fini sur l'anneau local $({\rm gr }^{\sharp} \, {\cal D}_{X  }[s_1, \ldots ,s_p]_{x_0})_{{\cal M}_{x_0,a,b}}$. Sa dimension est ind\'ependante de le bonne filtration di\`ese de $M$ et  co\"{\i}ncide avec  $ {\rm dim }_{x_0,a,b}    \,{\rm car}^{\sharp} \, M$ la dimension en $(x_0,a,b)$ de ${\rm car}^{\sharp} M$ 
(assertion analogue \`a la remarque 12 Chap\^{\i}tre  5,  \cite{G-M})~:
$$ {\rm dim } (({\rm gr }^{\sharp} \, M)_{x_0})_{{\cal M}_{x_0,a,b}}  = {\rm dim }_{x_0,a,b}    \,{\rm car}^{\sharp} \, M  \; .$$

 Comme les fibres de la restriction
de $\pi$ \`a ${\rm car}^{\sharp} M$ sont coniques~:
$$   {\rm dim }_{x_0,0,0}    \,{\rm car}^{\sharp}\,  M =  {\rm sup_{ (a,b) \in {\bf C}^n \times  {\bf C}^p  } } \; {\rm dim }_{x_0,a,b}    \,{\rm car}^{\sharp} \, M \; .$$
Le module   $({\rm gr }^{\sharp}\,  M)_{x_0}$ est gradu\'e et ${\cal M}_{x_0,0,0} $ est le seul id\'eal maximal gradu\'e de  ${\rm gr }^{\sharp}\, {\cal D}_{X  }[s_1, \ldots ,s_p]_{x_0}$. Il s'en suit que~:
$${\rm dim }  ({\rm gr }^{\sharp} \,  M)_{x_0} = {\rm dim } (({\rm gr }^{\sharp} \,ÊM)_{x_0})_{{\cal M}_{x_0,0,0}}  \; . $$

 Suivant J.-P Serre  (\cite{Se},
 Chap\^{\i}tre 4), un anneau $A$  commutatif noeth\'erien est r\'egulier si la borne sup\'erieure ${\rm gldh} \, (A)$ des entiers $k$ tels que $   Ext^k(M,N)\neq 0 $ pour un couple $M,N$ de $A$-module est finie. 
 Cette borne sup\'erieure est appel\'ee la  dimension homologique globale de $A$.
 L'anneau des s\'eries convergentes ${\bf C}\{x_1, \ldots ,x_n\} $ est r\'egulier de dimension homologique globale $n$. Des propri\'et\'es de transfert de la r\'egularit\'e (\cite{Se},
 Chap\^{\i}tre 4 proposition 25), il r\'esulte que l'anneau ${\bf C}\{ x_1, \ldots ,x_n \}[s_1 , \ldots ,s_p, \xi _1  \ldots \xi _n ] $ est r\'egulier de dimension 
 $2  n +p$. L'anneau localis\'e $({\rm gr }^{\sharp} \, {\cal D}_{X  }[s_1, \ldots ,s_p]_{x_0})_{{\cal M}_{x_0,0,0}}$ est alors un anneau local
 r\'egulier de dimension $2 \, {\rm dim } \, X +p$ (\cite{Se}, Chap\^{\i}tre 4 proposition 23). C'est donc un anneau de Cohen-Macaulay 
 (\cite{Se}, Chap\^{\i}tre 4  paragraphe D  corollaire 3): la longueur des suites r\'eguli\`eres maximales form\'ees d'\'el\'ements de son id\'eal maximal co\"{\i}ncide avec la dimension de l'anneau.
 Suivant D. Rees \cite{Re}, si $A$ est un anneau commutatif noeth\'erien et $E$ un A-module commutatif, nous appelons nombre grade de $E$ l'entier :
 $$ {\rm grade}\,  E  = {\rm inf}  \{ i \in {\bf N} \; ; \;   Ext^i_A( E,A ) \neq 0 \} \; .$$
Si $I =  {\rm ann} \,  E $ est l'annutateur de $E$, $ {\rm grade} \,  E =   {\rm grade} \, (A/I) $ et est le nombre maximum d'\'el\'ements d'une suite r\'eguli\`ere de
 $A$ form\'ee d'\'el\'ements de $I$  (\cite{Ma}, Chap\^{\i}tre 6, paragraphe 15 D).\\
 
 Si $M$ est un ${\cal D}_X[s_1, \ldots ,s_p] $-Module coh\'erent, pour tout $x_0 \in X$, les id\'eaux associ\'es \`a $({\rm gr }^{\sharp} M)_{x_0}$
 sont contenue dans ${\cal M}_{x_0,0,0}$. Nous obtenons alors en utilisant la platitude    et la commutation des Ext \`a la localisation pour les modules de type finie (\cite{Se}, Chap\^{\i}tre 4 proposition 18)~:
 $$ {\rm grade}\, ({\rm gr }^{\sharp} M)_{x_0}  =  {\rm grade}\, (({\rm gr }^{\sharp} M)_{x_0} )_{{\cal M}_{x_0,0,0}}  \; .$$
 Rappelons que si $E$ est un $A$-module de type fini sur un anneau $A$ local commutatif noeth\'erien Cohen-Macaulay (\cite{Ma}, Chap\^{\i}tre 6 theorem 31).~:
 $$  {\rm grade}\,  E + {\rm dim }\, \frac{A}{  {\rm ann} \,  E  }  =  {\rm dim }\, A \;    $$
 Il en r\'esulte~:
 
 \begin{proposition} \label{pgfd} Soit $M$ est un ${\cal D}_X[s_1, \ldots ,s_p] $-Module coh\'erent,  pour tout $x_0 \in X$ :
 $$ {\rm dim }\, ({\rm gr }^{\sharp} M)_{x_0} = 2 \, {\rm dim } \, X +p - {\rm grade}\, ({\rm gr }^{\sharp} M)_{x_0} = {\rm dim }_{x_0,0,0}    \,{\rm car}^{\sharp} M \; .$$
 \end{proposition}
 
 \subsection{Filtration relative
  de ${\cal D}_X[s_1, \ldots ,s_p] $}
  \label{ssfr}
La filtration $( {\cal D}^{\rm rel}_X[s_1, \ldots ,s_p](k) )_{k \in {\bf N}}$ de ${\cal D}_X[s_1, \ldots ,s_p] $ est 
d\'efinie comme suit : un op\'erateur $P \in {\cal D}^{\rm rel}_X[s_1, \ldots ,s_p](k) $ s'il s'\'ecrit localement~:
$$  P = \sum_{   \mid \beta \mid \leq k}  a_{ \beta}(x,s)  \partial ^{\beta} \quad {\rm  o\grave{u} } 
\quad  a_{ \beta}\in {\cal O}_X[s_1, \ldots ,s_p] 
 \;   {\rm  et}   \; \partial ^{\beta} = \frac{ \partial ^{\beta _1} }{ \partial x_1 }  \cdots \frac{ \partial ^{\beta _n} }{\partial x_n}  \; .$$
 Nous appelons cette filtration   la
   \underline{filtration relative} de ${\cal D}_X[s_1, \ldots ,s_p] $. \\
   
Nous d\'esignons par ${\rm gr }^{\rm rel} \, {\cal D}_X[s_1, \ldots ,s_p]$ le gradu\'e de cette filtration. 
Ce gradu\'e  s'identifie au sous-faisceau 
de $\pi _{\ast} ( {\cal O}_{T^{\ast}X \times {\bf C}^p} )$ des fonctions analytiques sur $T^{\ast} X \times {\bf C}^p$ polynomiales par rapport aux fibres de $\pi$. Il
est isomorphe localement   au faisceau d'anneaux commutatifs  gradu\'es $  {\cal O}_ X [s_1, \ldots ,s_p, \xi _1,  \ldots ,\xi _n ] $.\\

Un op\'erateur $P \in{\cal D}^{\rm rel}_X [s_1, \ldots ,s_p](k)  - {\cal D}^{\rm rel}_X[s_1, \ldots ,s_p](k-1)$   modulo $ {\cal D}^{\rm rel}_X[s_1, \ldots ,s_p](k-1)$  d\'efinit une section de ${\rm gr }^{\rm rel} {\cal D}_X[s_1, \ldots ,s_p]$
appel\'ee symbole relatif   de $P$ et not\'e $\sigma ^{\rm rel}(P)$.  Localement~: \\
$$ {\rm Si  } 
\quad  P  =  \sum_{   \mid \beta \mid \leq k}  a_{{\alpha},\beta}(x,s)  \partial ^{\beta} \quad , \quad 
 \sigma ^{\rm rel}(P) (x,\xi,s) = \sum_{   \mid \beta \mid = k}  a_{ \beta}(x,s)  \xi ^{\beta} \; .$$
Le faisceau d'anneaux ${\rm gr }^{\rm rel} \, {\cal D}_X[s_1, \ldots ,s_p]$ s'identifie \`a  $  {\cal O}_ X [s_1 , \ldots ,s_p, \xi _1  \ldots \xi _n ] $. C'est un faisceau coh\'erent d'anneaux commutatifs.  
 Soit $M$  un  ${\cal D}_X[s_1, \ldots ,s_p] $-Module coh\'erent et $(M_k)_{k \in {\bf Z} }$ une filtration de $M$ pour la filtration relative. Nous disons que cette filtration est bonne si
 localement, il existe des sections $m_1, \ldots ,m_l$ de $M$  et des entiers relatifs $k_1, \ldots , k_p$ tels que~:
  $$M_k = \sum_{i=1}^p {\cal D}^{\rm rel}_X[s_1, \ldots ,s_p](k-k_i)   m_i \; .$$
 Nous notons ${\rm gr }^{\rm rel} \, M$ le gradu\'e de $M$ pour cette filtration et ${\rm Ann}\, {\rm gr }^{\rm rel} \, M$ son annulateur comme 
${\rm gr }^{\rm rel} \, {\cal D}_X[s_1, \ldots ,s_p]$-Module. La racine de ${\rm Ann}\, {\rm gr }^{\rm rel} \, M$ est un id\'eal coh\'erent
de ${\rm gr }^{\rm rel} \, {\cal D}_X[s_1, \ldots ,s_p]$ ind\'ependant  des bonnes filtrations relatives de $M$. Nous
notons $J^{\rm rel} (M)$ cet id\'eal. Nous appelons \underline{vari\'et\'e  caract\'eristique  relative} de $M$
le  sous-espace analytique de $T^{\ast}X \times {\bf C}^p  $ d\'efini par  $J^{\rm rel} (M)$ et le
  notons   ${\rm car}^{\rm rel} \, M$.\\

Soit $p : X \times {\bf C}^p \rightarrow X$ la  projection sur $X$. 
Pour tout ${\cal D}_X[s_1, \ldots ,s_p] $-Module coh\'erent $M$, posons~:
$ M^{\rm an} = {\cal O}_{ X \times {\bf C}^p} \otimes_{p^{-1} {\cal O}_X[s_1, \ldots ,s_p]} p^{-1} M$. Le faisceau
${\cal D}_X[s_1, \ldots ,s_p]^{\rm an}  =   {\cal O}_{ X \times {\bf C}^p } \otimes_{ p^{-1} {\cal O}_X[s_1, \ldots ,s_p] } p^{-1} {\cal D}_X[s_1, \ldots ,s_p] $
s'identifie au faisceau ${\cal D}_{ X \times {\bf C}^p /  {\bf C}^p }$ des op\'erateurs diff\'erentiels relatifs. La filtration 
naturelle
de ${\cal D}_{ X \times {\bf C}^p /  {\bf C}^p }$ par l'ordre des d\'erivations n'est autre que :
$$ {\cal O}_{ X \times {\bf C}^p } \otimes_{ p^{-1} {\cal O}_X[s_1, \ldots ,s_p] } p^{-1}  {\cal D}^{\rm rel}_X [s_1, \ldots ,s_p](k) \; .$$
 Le faisceau $ M^{\rm an}$ est   un  ${\cal D}_{ X \times {\bf C}^p /  {\bf C}^p }$-Module
coh\'erent. Si $(M_k)_{k \in {\bf Z} }$  est une bonne filtration relative de $M$, $ {\cal O}_{ X \times {\bf C}^p } \otimes_{ p^{-1} {\cal O}_X[s_1, \ldots ,s_p] } p^{-1} M_k$
est une bonne filtration de $ M^{\rm an}$. Le gradu\'e de cette filtration s'identifie \`a 
$ {\cal O}_{ X \times {\bf C}^p } \otimes_{ p^{-1} {\cal O}_X[s_1, \ldots ,s_p] } p^{-1} {\rm gr }^{\rm rel} M $. Par platitude de  $ {\cal O}_{ X \times {\bf C}^p }$
sur $p^{-1} {\cal O}_X[s_1, \ldots ,s_p]  $, son annulateur  est $ {\cal O}_{ X \times {\bf C}^p } \otimes_{ p^{-1} {\cal O}_X[s_1, \ldots ,s_p] } p^{-1} {\rm Ann}\, {\rm gr }^{\rm rel} M $ 
(voir \cite{Bo.AC}   Chap\^{\i}tre 1, paragraphe 2, corollaire 2). Nous notons $J (M^{\rm an} ) $ la racine de cet annulateur qui ne d\'epend pas des bonnes filtrations de $M^{\rm an}$.
Le sous-espace analytique  ${\rm car}  \, M^{\rm an}$   de  $T^{\ast}X \times {\bf C}^p  $ d\'efini par $J (M^{\rm an} ) $  co\"{\i}ncide avec ${\rm car}^{\rm rel} \, M$.\\

Les faisceaux de Modules $ {\rm gr }^{\rm rel} \, M$ et $ {\rm gr }^{\rm rel} \, M^{\rm an}$ sur ${\rm gr }^{\rm rel} \, {\cal D}_X[s_1, \ldots ,s_p]$ et
$ {\rm gr } \, {\cal D}_{ X \times {\bf C}^p /  {\bf C}^p } $ sont coh\'erents (\cite{Bj2} Appendice III, theorem 2.17). Par platitude, pour tout entier $i$~:
$$ Ext^i_{ {\rm gr } \, {\cal D}_{ \frac{X \times {\bf C}^p }{ {\bf C}^p   }}}    ( {\rm gr }  \,  M^{\rm an} ,  {\rm gr } \, {\cal D}_{ \frac{X \times {\bf C}^p }{ {\bf C}^p  } } )  \quad et \quad Ext^i_{ {\rm gr }^{\rm rel} \, {\cal D}_X[s_1, \ldots ,s_p]} ( {\rm gr }^{\rm rel} \, M,{\rm gr }^{\rm rel} \, {\cal D}_X[s_1, \ldots ,s_p] )^{\rm an}   $$ 
  sont isomorphes.
Il en r\'esulte que pour tout $x_0 \in X$, il existe $U$ voisinage de $x_0$ tel que pour tout $(x,b) \in U \times {\bf C}^p$~:
$$ {\rm grade}\, ( {\rm gr }   \,M^{\rm an}  )_{x,b} \geq {\rm grade}\, ( {\rm gr }^{\rm rel} \, M )_{x_0}  $$
La m\^eme preuve que dans le cas de la filtration di\`ese  portant sur le fait qu'un anneau  local r\'egulier   est Cohen-Macaulay montre :
$$ \begin{array}{rcl} {\rm grade}\, ( {\rm gr }  \, M^{\rm an}  )_{x,b}  &=&  2 \, {\rm dim } \, X +p - {\rm dim }_{x,0,b}    \,{\rm car}^{\rm rel} \, M \\
 &= &2 \, {\rm dim } \, X +p - {\rm sup}_{ a \in {\bf C}^n} {\rm dim }_{x,a,b}    \,{\rm car}^{\rm rel} \,M \; . \end{array}  $$  
 Nous obtenons donc :
\begin{remarque}  \label{ridg}Si $M$ est un un  ${\cal D}_X[s_1, \ldots ,s_p] $-Module coh\'erent, pour tout $x_0 \in X$, il existe $U$ voisinage de $x_0$
tel que pour tout $(x,b) \in U \times {\bf C}^p$ :
$$ \begin{array}{rcl} {\rm dim }_{x,0,b}    \,{\rm car}^{\rm rel} \, M   &=&  {\rm sup}_{ a \in {\bf C}^n} \, {\rm dim }_{x,a,b}    \,{\rm car}^{\rm rel} \, M \\
 & \leq  &2 \, {\rm dim } \, X +p -   {\rm grade}\, ( {\rm gr }^{\rm rel} \, M )_{x_0}  \; . \end{array}  $$ 
\end{remarque} 

Une  difficult\'e de  l'anneau  $ {\rm gr }^{\rm rel} {\cal D}_X[s_1, \ldots ,s_p] $ est d'une part que ses fibres 
  n'ont pas m\^eme hauteur et d'autre part
le  d\'efaut de fid\'ele platitude de
   ${\cal O}_{ X \times {\bf C}^p }$  sur  $p^{-1} {\cal O}_X[s_1, \ldots ,s_p]$. Par exemple, l'anneau quotient~:
   $$  L = \frac{{\bf C}\{ x_1, \cdots ,x_n\}[s_1, \ldots ,s_p]  }{(1-s_1x_1)} \neq 0 \; , $$
mais~:
$ {\bf C} \{ x_1, \cdots ,x_n,s_1 ,\ldots, s_p \}  \otimes_{{\bf C}\{ x_1, \cdots ,x_n\}[s_1, \ldots ,s_p]  }  L =0 $. 
Cette difficult\'e \'etait  d\'eja soulign\'ee dans \cite{M-N}. \\

\begin{lemme} (Th\'eor\`eme des z\'eros relatifs) si $N$ est un $ {\rm gr }^{\rm rel} \, {\cal D}_X[s_1, \ldots ,s_p]$-Module coh\'erent, pour tout $x_0$ in $X$~:
$$ N_{x_0} = 0 \;  \Longleftrightarrow    \;  \exists \, U \;   {\rm voisinage \; de}\;   x_0   \;   {\rm tel \; que }\;  N^{\rm an}\, _{\mid p^{-1} (U)} = 0 \; .$$
\end{lemme}

\noindent Preuve : $\Rightarrow$ est facile. Nous montrons l'autre implication dans le sous-paragraphe \ref{thzr}.

\begin{proposition} \label{pgfr} Soit $M$  un  ${\cal D}_X[s_1, \ldots ,s_p] $-Module coh\'erent et $x_0 \in X$~:
 $$ \begin{array}{rcl} {\rm grade} \, ( {\rm gr }^{\rm rel}  \,  M   )_{x_0}  &=& {\sup}_{ U \, {\rm voisinage \; de } \, x_0 }  {\inf}_{  (x,b) \in p^{-1}(U)    }  \, 
 {\rm grade}\, ( {\rm gr }   M^{\rm an}  )_{x,a}
 \\ &= &2 \, {\rm dim } \, X +p -  {\inf}_{U \, {\rm voisinage \; de } \, x_0} \, {\sup}_{ (x,b) \in p^{-1}(U)   }  {\rm dim }_{x_0,0,b}    \,{\rm car}^{\rm rel} \, M \; . 
 \end{array}      $$
 \end{proposition}
 
 \noindent Preuve : Soit $x_0 \in X$ et $U$ voisinage de $x_0$. Posons~:
 $$ k = {\inf}_{  (x,b) \in p^{-1}(U)    }  \, 
 {\rm grade}\, ( {\rm gr }   M^{\rm an}  )_{x,b}
 \; .$$
 Pour $i<k$, pour tout $(x,b) \in p^{-1}(U) $~:
 $$ Ext^i_{ {\rm gr } \, {\cal D}_{ \frac{X \times {\bf C}^p }{ {\bf C}^p   }}}    ( {\rm gr }  \,  M^{\rm an} ,  {\rm gr } \, {\cal D}_{ \frac{X \times {\bf C}^p }{ {\bf C}^p  } } ) _{x,b}=0 \; .$$
 Il r\'esulte du lemme que pour $i<k$~:
 $$Ext^i_{ {\rm gr }^{\rm rel} \, {\cal D}_X[s_1, \ldots ,s_p]} ( {\rm gr }^{\rm rel} \, M,{\rm gr }^{\rm rel} \, {\cal D}_X[s_1, \ldots ,s_p] )_{x_0} = 0 \; .$$
 Donc, pour tout $U$ voisinage de $x_0$~: 
 $${\rm grade} \, ( {\rm gr }^{\rm rel}  \,  M   )_{x_0}  \geq  {\inf}_{  (x,b) \in p^{-1}(U )   }  \, 
 {\rm grade}\, ( {\rm gr }   M^{\rm an}  )_{x,b} \; .$$
 et donc~: 
 $${\rm grade} \, ( {\rm gr }^{\rm rel}  \,  M   )_{x_0}  \geq  {\sup}_{ U \, {\rm voisinage \; de } \, x_0 }   {\inf}_{  (x,b) \in p^{-1}(U)    }  \, 
 {\rm grade}\, ( {\rm gr }   M^{\rm an}  )_{x,b} \; .$$

 D'autre part, voir preuve de la  remarque \ref{ridg}, il existe $U_0$ voisinage de $x_0$ tel que pour tout $(x,a) \in U_0 \times {\bf C}^p$~:
$$  {\rm grade}\, ( {\rm gr }^{\rm rel} \, M )_{x_0}     \leq  {\rm grade}\, ( {\rm gr }   \,M^{\rm an}  )_{x,b} $$
 
Ainsi, sur ce voisinage $U_0$~:
$$ {\rm grade}\, ( {\rm gr }^{\rm rel} \, M )_{x_0}  \leq    {\inf}_{  (x,b) \in p^{-1}(U_0)    }  \, 
 {\rm grade}\, ( {\rm gr }   M^{\rm an}  )_{x,b}
 \; ,$$
 ce qui montre la proposition.

 \subsection{Dimension d'un ${\cal D}_X[s_1, \ldots ,s_p] $-Module coh\'erent}
 
 Soit $B$ un anneau non n\'ecessairement commutatif, nous disons qu'il est r\'egulier s'il est noeth\'erien \`a gauche et \`a droite et de dimension homologique
 globale finie. Si $E$ est un $B$-module de type fini, par exemple  \`a gauche, nous appelons   nombre  grade de $E$, l'entier~:
 $$ {\rm grade}\,  E  = {\rm inf}  \{ i \in {\bf N} \; ; \;   Ext^i_B( E,B ) \neq 0 \} \; .$$

 Suivant J.-E. Bj\"{o}rk (\cite{Bj2}  Appendice IV), nous disons qu'un anneau r\'egulier $B$ v\'erifie la condition d'Auslander si pour tout 
  $B$-module \`a gauche $E$ de type fini   et tout sous-module $N$ \`a  droite de $Ext^k_B(E,B)$, le grade de $N$ est sup\'erieure ou \'egal \`a $k$.\\
  
 Tout anneau $A$ r\'egulier commutatif v\'erifie la condition d'Auslander (\cite{A-B}). Rappelons une preuve de la d\'emonstration.  Tout d'abord, nous nous ramenons
 facilement  au cas o\`u $A$  est local.
 Dans ce cas, pour tout module $N$ de type fini, $ {\rm grade}\, N=  {\rm dim } \, A -  {\rm dim } \,  N$. En utilisant cette \'egalit\'e, il suffit de   
   montrer  que 
  si $M$ est un module de type fini,  $ {\rm grade}\, Ext^k_A(M,A)\geq k$. Soit $q$ un id\'eal premier appartenant au support de $Ext^k_A(M,A)$. Alors,
  $ Ext^k_{A_q}(M_q,A_q)=Ext^k_A(M,A)_q\neq 0$. Ainsi, $gl {\rm dim } \, A_q$  la dimension homologique globale de $A_q$ est sup\'erieure ou \'egale \`a $k$. Or l'anneau
  $A_q$ est r\'egulier. Donc, suivant \cite{Se}, Chap\^{\i}tre IV, corollaire 2, $gl {\rm dim } \, A_q = {\rm dim } \, A_q$.    Choissisons $q$ tel que ${\rm dim } \, A/q = {\rm dim } \,Ext^k_A(M,A)$.
  Comme $A$ est Cohen Macaulay ${\rm ht}\,  (q) + {\rm dim } \, A/q = {\rm dim } \, A$ o\`u ${\rm ht}\,  (q) $ d\'esigne la hauteur de l'id\'eal, $q$. Nous obtenons~:
  $${\rm grade}\, Ext^k_A(M,A) = {\rm dim } \, A - {\rm dim } \, \frac{A}{q} ={\rm ht}\,  (q)  =  {\rm dim } \, A_q = gl {\rm dim } \, A_q \geq k \; .$$
  
  Cette propri\'et\'e d'Auslander est \`a la base des propri\'et\'es du nombre grade et de la filtration d'un module par son complexe bidualisant.\\
  
  Suivant Bj\"ork (\cite{Bj2}, appendice IV, theorem 4.15), si $B$ est un anneau filtr\'e positivement dont le gradu\'e est r\'egulier et v\'erifie la condition d'Auslander, alors
  $B$ est r\'egulier et v\'erifie la condition  d'Auslander. De plus, pour tout $B$ module de type fini \`a gauche~:
  $$ {\rm grade}_B \,M =   {\rm grade}_{{\rm gr } \, B}  \, {\rm gr }   \,M \; ,$$  
  o\`u ${\rm gr }   \,M$ 
  d\'esigne le gradu\'e d'une bonne filtration de $M$.  L'application de ces r\'esultats donnent~:
  
  \begin{proposition} \label{pegr} Les fibres de ${\cal D}_X[s_1, \ldots ,s_p]  $ sont des anneaux r\'eguliers qui v\'erifient la condition d'Auslander. De plus, pour tout 
  ${\cal D}_X[s_1, \ldots ,s_p] $-Module coh\'erent  $M$ et $x_0 \in X$~:
  $${\rm grade  }  \, M_{x_0}   = {\rm grade  }  \,  ({\rm gr }^{\sharp} \, M)_{x_0}   = {\rm grade  }\,    ( {\rm gr }^{\rm rel} \, M)_{x_0}    \quad . $$
 \end{proposition}
  
  Compte-tenu des   propositions \ref{pgfd} et   \ref{pgfr}, nous obtenons le th\'eor\`eme suivant qui \'etabli le lien entre le grade d'un 
  ${\cal D}_X[s_1, \ldots ,s_p] $-Module coh\'erent et les dimensions de ses vari\'et\'es caract\'eristiques di\`ese et relative.
\begin{theoreme} Pour tout 
  ${\cal D}_X[s_1, \ldots ,s_p] $-Module coh\'erent  $M$ et $x_0 \in X$~:
$$ \begin{array}{rcl}  {\rm dim }_{x_0,0,0}    \,{\rm car}^{\sharp} M &=&{\inf}_{U \, {\rm voisinage \; de } \, x_0} \, {\sup}_{ (x,a) \in p^{-1}(U)   }  {\rm dim }_{x_0,a,0}    \,{\rm car}^{\rm rel} \, M \\ &= & 2 \, {\rm dim } \, X +p -  {\rm grade }  \, M_{x_0}  \; . 
 \end{array}      $$
 \end{theoreme}

  \subsection{Th\'eor\`eme des z\'eros relatifs}
  
  \begin{proposition}\label{thzr} Soit $X$ une vari\'et\'e analytique et $I$ un id\'eal coh\'erent de  ${\cal O}_X[s_1, \ldots ,s_p] $ tel que~: 
  $$ {\rm V} (I) = \{ (x,s) \in X \times {\bf C}^p    \; ; \; \forall \, g \in I \; : \; g(x,s) = 0 \}  =  \emptyset \; ,  $$
  alors, $I = {\cal O}_X[s_1, \ldots ,s_p] $.
  \end{proposition}
 
 \noindent {\bf Preuve :  } La question est locale sur $X$. Nous pouvons supposer que $X$ est un voisinage de l'origine dans $ {\bf C}^n$
 et $I = (g_1, \ldots ,g_l) $ o\`u $g_i \in {\cal O}_X(X)[s_1, \ldots ,s_p] $. Par hypoth\`ese   le syst\`eme d'\'equations~:
 $$ g_1(x,s) = \cdots = g_l(x,s)=0 \quad , \quad (x,s) \in  X \times {\bf C}^p $$ 
 n'a pas de solution. Nous devons montrer qu'il existe $U$ voisinage de l'origine dans $ {\bf C}^n$ et $a_1, \ldots ,a_l \in  {\cal O}_X(U)[s_1, \ldots ,s_p]$ tels que~:
 $$ 1 = \sum_{i=1}^{l}  a_i\,  {g_i}_{\mid U} \; .$$
Soit un entier $d$ tel que pour tout $x \in X$, les degr\'es en $s$ des polyn\^omes $g_i(x,s)$ soient  born\'es par $d$.
Chaque $g_i$ s'\'ecrit $\sum_{\mid \beta \mid \leq d }  g_{i, \beta}(x)s^{\beta}$ o\`u $g_{i, \beta} \in {\cal O}_X$.  D'apr\`es le th\'eor\`eme des z\'eros effectif   (\cite{H},\cite{B}), il existe un entier $N$ tel  que pour tout  $x\in X$, il existe des polyn\^omes $a_1(x)(s),\ldots , a_l(x)(s)$ en $s$ de degr\'e  inf\'erieur  ou \'egal \` a $N$ tels que~: 
 $$ 1 = \sum_{i=1}^{l}  a_i (x)(s)\,  {g_i} (x,s) \; .$$

Consid\'erons alors l'\'equation  polynomiale    dont les inconnues sont  les   $a_{i,\alpha}$ o\`u $i \in \{ 1, \ldots , l\}$, $ \alpha \in {\bf N}^p$ et
$ \mid \alpha \mid \leq N$~:
$$ \sum_{i=1}^l  (  \sum_{\mid \alpha \mid \leq N } a_{i,\alpha} s^{\alpha }  ) ( \sum_{\mid \beta \mid \leq d }  g_{i, \beta}(x)s^{\beta}  ) = 1 \; , $$
Cette  \'equation   s'\'ecrit~:
$$\sum_{\mid \gamma \mid \leq N+d }  ( \sum_{i=1}^l  ( \sum^{ \alpha + \beta =  \gamma }_{ \mid \alpha \mid \leq N , \mid \beta \mid \leq d}   a_{i,\alpha} g_{i, \beta}(x) )  s^{\gamma }  = 1\; . $$ 
 Elle est \'equivalente au  syst\`eme  lin\'eaire   d'inconnues    les   $a_{i,\alpha}$  :
 $$\left[ \begin{array}{rcl}   \sum_{i=1}^l     g_{i,0}(x) a_{i,0} = 1 \\  \sum_{i=1}^l   \sum^{ \alpha + \beta =  \gamma\neq 0 }_{ \mid \alpha \mid \leq N , \mid \beta \mid \leq d}    g_{i, \beta}(x)   a_{i,\alpha}   = 0 \; .\end{array}   \right.$$
 Nous avons vu  que pour tout $x \in U$, ce  syst\`eme  a  des solutions.  Apr\`es avoir ordonn\'e les $\alpha$ et $\gamma$, ce  syst\`eme  s'\'ecrit~:
 $$ G(x) A = \left(  \begin{array}{c}    1\\ 0 \\ \vdots \\ 0 \end{array}    \right) \; .$$
 o\`u $A$ est une matrice colonne de coefficients $a_{i,\alpha}$ et $G(x)$ une matrice dont les coefficients non nuls sont des $g_{i,\beta}(x)$.\\
 
 Soit $\Delta _1$ un mineur de la matrice $G(x)$ de taille maximum $r_1$ non identiquement nul dans un voisinage $U$ de l'origine. Suivant les formules
 de Cramer, il existe des $b_{i,\alpha,1}(x)$ mineurs de taille $r_1 -1$ de $G(x)$ tels que pour tout $x \in U$ v\'erifiant  $ \Delta _1 (x) \neq 0 $~:
 $$  \sum_{i=1}^l    (  \sum_{\mid \alpha \mid \leq N }    \frac{b_{i,\alpha,1}(x)}{ \Delta _1 (x)  } s^{\alpha })   {g_i} (x,s)  =  1\; .$$ 
 Ainsi, pour tout $x\in U$ tels que  $ \Delta _1 (x) \neq 0 $~:
 $$ \Delta _1 (x) =  \sum_{i=1}^l    (  \sum_{\mid \alpha \mid \leq N }    b_{i,\alpha,1}(x) s^{\alpha })   {g_i} (x,s) \; .$$
 Dons, pour tout $x \in U$ :
  $$ \Delta _1^2  (x) =  \sum_{i=1}^l    (  \sum_{\mid \alpha \mid \leq N }    b_{i,\alpha,1}(x) \Delta _1 (x)s^{\alpha })   {g_i} (x,s) \; .$$
 Si $\Delta _1 (0) \neq 0  $, quitte \`a diminuer $U$, nous obtenons l'existence de $a_1; \ldots , a_p \in  {\cal O}_X(U)[s_1, \ldots ,s_p]$ tels que~:
 $$ 1 = \sum_{i=1}^{l}  a_i\,  {g_i}_{\mid U} \; .$$
Si $\Delta _1 (0)= 0  $ : Consid\'erons sur $\Delta _1 ^{-1} (0) \subset X $, un mineur $\Delta _2$  de la matrice $G(x)$ de taille maximum $r_2$
non identiquement nul dans un voisinage $U$ de l'origine.  Suivant les formules
 de Cramer, il existe des $b_{i,\alpha,2}(x)$ mineurs de taille $r_2 -1$ de $G(x)$ tels que pour tout $x \in U$ v\'erifiant   $ \Delta _1 (x) = 0 $ et $ \Delta _2 (x) \neq 0 $~:
  $$\sum_{i=1}^l    (  \sum_{\mid \alpha \mid \leq N }    \frac{b_{i,\alpha,2}(x)}{ \Delta _2 (x)  } s^{\alpha })   {g_i} (x,s)  =  1\; .  $$ 
  Il en r\'esulte :
  $$ \Delta _2^2(x) - \sum_{i=1}^l    (  \sum_{\mid \alpha \mid \leq N }     b_{i,\alpha,2}(x)  \Delta _2 (x)    s^{\alpha })   {g_i} (x,s) $$ 
  nul sur $\Delta _1 ^{-1} (0)$.  D'apr\'es le th\'eor\`eme des z\'eros analytiques,  il existe quitte \`a diminuer $U$  un entier $k$
   tel que ~: 
   $$ (\Delta _2^2(x) - \sum_{i=1}^l    (  \sum_{\mid \alpha \mid \leq N }     b_{i,\alpha,2}(x)  \Delta _2 (x)    s^{\alpha })   {g_i} (x,s))^k 
   \in {\cal O}_X(U)[s_1, \ldots ,s_p] ( \Delta _1^2  (x) )\; .$$
   Vu l'\'ecriture de $\Delta _1^2$, nous en d\'eduisons quitte \`a diminuer $U$,  l'existence de $a_1, \ldots , a_p \in  {\cal O}_X(U)[s_1, \ldots ,s_p]$ tels que~:
$$\Delta _2^{2k}(x) = \sum_{i=1}^{l}  a_i\,  {g_i}_{\mid U} \; .$$
Si $\Delta _2 (0) \neq 0  $, nous terminons comme pr\'ecedemment. Sinon, nous it\'erons. Le processus aboutit car les mineurs de taille $1$ de $G(x)$ 
ne sont pas tous nuls en z\'ero.

   \begin{corollaire} Soit $X$ une vari\'et\'e analytique complexe,  $N$ un  ${\cal O}_X[s_1, \ldots ,s_p] $-Module coh\'erent tel que
   $N^{\rm an}= {\cal O}_{X \times {\bf C}^p }\otimes_{p^{-1}  {\cal O}_X[s_1, \ldots ,s_p]} p^{-1} \, N $ soit nul o\`u $p$ d\'esigne la projection de  
   $X \times {\bf C}^p$ sur ${\bf C}^p$. Alors, $N=0$.   
      \end{corollaire}
   
  \noindent {\bf Preuve :  }  $({\rm Ann}\,  N)^{\rm an}  =  {\rm Ann}\, ( N^{\rm an} )$.  Il suffit alors d'appliquer la proposition \ref{thzr} \`a l'id\'eal
  ${\rm Ann}\,  N $.
   
      \begin{corollaire} \label{ctzr} Soit $X$ une vari\'et\'e analytique et $I$ un id\'eal coh\'erent de  ${\cal O}_X[s_1, \ldots ,s_p] $. Soit $h \in  {\cal O}_X[s_1, \ldots ,s_p] $
      tel que $h$ s'annule sur 
  $$ {\rm V} (I) = \{ (x,s) \in X \times {\bf C}^p    \; ; \; \forall \, g \in I \; : \; g(x,s) = 0 \}   \; .  $$
  Alors, localement au voisinage de tout point de $X$, il existe $k$ tel que $h^k \in I$.
\end{corollaire}

  \noindent {\bf Preuve :}  L'astuce de Rabinowitch s'adapte sans probl\`eme.

   \newpage   
      
   \section{${\cal D}_X[s_1, \ldots ,s_p] $-Module major\'e par une lagrangienne}
   \label{s2}
     
   \subsection{Cons\'equence du th\'eor\`eme d'involutivit\'e} \label{scti}
   
   Consid\'erons tout d'abord sur ${\cal D}_X[s_1, \ldots ,s_p] $ la filtration di\`ese. Le faisceau ${\cal D}_X[s_1, \ldots ,s_p] $
   est un faisceau de ${\bf Q}$-alg\`ebres et ${\rm gr }^{\sharp} \, {\cal D}_X[s_1, \ldots ,s_p]$ est un faisceau d'anneaux commutatifs. Si $P \in
   {\cal D}_X^{\sharp}[s_1, \ldots ,s_p](k)$ et $Q \in      {\cal D}_X^{\sharp}[s_1, \ldots ,s_p](l) $ : $PQ-QP \in {\cal D}_X^{\sharp}[s_1, \ldots ,s_p](k+l-1)$.
   Dans un syst\`eme de coordonn\'des locales si $PQ-QP \notin {\cal D}_X^{\sharp}[s_1, \ldots ,s_p](k+l-2)$~:
   $$ \sigma^{\sharp}   (PQ-QP) = \sum_{i=1}^n \frac{\partial \sigma ( P) }{ \partial \xi _i} \frac{\partial \sigma ( Q) }{ \partial x_i} - 
   \frac{\partial \sigma ( P) }{ \partial \xi i}\frac{\partial \sigma ( Q) }{ \partial x_i }
   \; .$$
   Cette formule \'etend le crochet de Poisson d\'efini sur les symboles des op\'erateurs de ${\cal D}_X$. Si $\alpha$ et $\beta$
   sont des fonctions sur ${\rm gr }^{\sharp} \, {\cal D}_X[s_1, \ldots ,s_p]$, nous appelons crochet de Poisson de $\alpha$ et $\beta$~:
   $$ \{  \alpha, \beta \} = \sum_{i=1}^n \frac{\partial \alpha }{ \partial \xi _i} \frac{\partial \beta }{ \partial x_i} - 
   \frac{\partial \alpha }{ \partial \xi i}\frac{\partial \beta}{ \partial x_i }
   \; .$$
   Le th\'eor\`eme de O. Gabber \cite{Ga} assure que si $M$ est un  ${\cal D}_X[s_1, \ldots ,s_p] $-Module coh\'erent \`a gauche,
   $J^{\sharp} (M)$ la racine de l'annulateur du gradu\'e de $M$ pour toutes  bonnes filtrations di\`eses est stable par crochet de Poisson.\\

   De m\^eme, si nous consid\'erons sur    ${\cal D}_X[s_1, \ldots ,s_p] $ la filtration  relative,    ${\rm gr }^{\rm rel} \, {\cal D}_X[s_1, \ldots ,s_p]$ est un
   encore un  faisceau d'anneaux commutatifs. Par le m\^eme th\'eor\`eme, nous obtenons que pour tout 
   ${\cal D}_X[s_1, \ldots ,s_p] $-Module coh\'erent \`a gauche $M$,
   $J^{\rm rel } (M)$ la racine de l'annulateur du gradu\'e de $M$ pour toutes  bonnes filtrations relatives est stable par crochet de Poisson.\\
   
  Soit $\pi _2 : T^{\ast }X \times {\bf C}^p \rightarrow  {\bf C}^p$, la projection sur   $ {\bf C}^p$. Pour tout sous-ensemble $Z$ de $T^{\ast }X \times {\bf C}^p$
   et $c \in {\bf C}^p$, nous notons $Z(c)$ la fibre au dessus de $c$ de la restriction de $\pi _2$ \`a $Z$. 
   
   \begin{proposition} Soit $M$ un ${\cal D}_X[s_1, \ldots ,s_p] $-Module coh\'erent. Pour tout $c \in  {\bf C}^p $, les fibres non vides
   $({\rm car }^{\sharp}\, M) (c)$ et    $({\rm car }^{\rm rel}\, M) (c)$   ont leurs composantes irr\'eductibes de dimension au moins \'egales \`a la dimension de $X$.
   \end{proposition}
   
     \noindent {\bf Preuve :  } En un point lisse de     ${\rm car }^{\sharp}\, M $ o\`u la restriction de $\pi _2$ est de rang localement constant, la fibre     $({\rm car }^{\sharp}\, M) (c)$ 
     est lisse r\'eduite. Par le th\'eor\`eme d'involutivit\'e,    $({\rm car }^{\sharp}\, M) (c)$  que nous identifions  \` a un sous-espace de $T^{\ast }X$ est involutive. Elle est donc de dimension sup\'erieure ou \'egale
     \`a la dimension de $X$. Par semi-continuit\'e  de la dimension des fibres, nous en d\'eduisons la propri\'et\'e annonc\'ee  sur $({\rm car }^{\sharp}\, M) (c)$.
     Le m\^eme raisonnement s'applique pour $( {\rm car }^{\rm rel}\, M ) (c)$.
     
     \subsection{${\cal D}_X[s_1, \ldots ,s_p] $-Module major\'e  par une lagrangienne}

   \begin{proposition} \label{pmvc} Pour tout    ${\cal D}_X[s_1, \ldots ,s_p] $-Module coh\'erent $M$~:
   $${\rm car }^{\rm rel}\, M \subset  ({\rm car }^{\sharp}\, M )(0) \times {\bf C}^p  \; . $$ 
   \end{proposition}

\noindent {\bf Preuve :  } Observons que si $P \in {\cal D}_X[s_1, \ldots ,s_p] $, $\sigma ^{\sharp}(P)(x,\xi,0)\neq 0$
implique $\sigma ^{\sharp}(P)(x,\xi,0) = \sigma ^{\rm rel}(P)(x,\xi,s)$. Pour d\'emontrer la proposition, nous pouvons supposer 
que $M$ est le quotient de ${\cal D}_X[s_1, \ldots ,s_p] $ par un id\'eal \`a gauche coh\'erent $I$. Soit 
Nous d\'eduisons de notre observation~:
$$ {\rm car }^{\rm rel}\, M \subset  V (\{ \sigma ^{\sharp}(P)(x,\xi,0) \; ; \; P \in I  \}  )=  ({\rm car }^{\sharp}\, M )(0) \times {\bf C}^p\; . $$

     \begin{definition}      Soit $M$ un  ${\cal D}_X[s_1, \ldots ,s_p] $-Module coh\'erent. Nous disons que
   $M$ est major\'e  par une lagrangienne $ \Lambda $   si $  {\rm car }^{\rm rel}\, M \subset \Lambda \times  {\bf C}^p $
   o\`u   $ \Lambda $ est une    sous-vari\'et\'e
   lagrangienne conique \'eventuellement singuli\`ere de $T^{\ast} X$.      \end{definition}  
       
       Suivant la proposition  \ref{pmvc}~:
       
       \begin{remarque} Si $({\rm car}^{\sharp} \, M )(0)$  est une vari\'et\'e
     lagrangienne de $T^{\ast} X$, alors    $M$ est major\'e  par la lagrangienne    $({\rm car}^{\sharp} \, M )(0)$   \end{remarque}
       
       Si $M$ est un ${\cal D}_X[s_1, \ldots ,s_p] $-Module coh\'erent, pour tout $c \in {\bf C}^p$, notons~:
       $$ M(c) = \frac{M}{\sum_{i=1}^p (s_i - c_i)M}\; .$$
     Ces ensembles  $M(c)$ sont naturellement des ${\cal D}_X$-Modules coh\'erents et nous    notons ${ \rm car}_{ {\cal D}_X  }\,  M(c) $ leurs vari\'et\'es caract\'eristiques
     comme ${\cal D}_X$-Modules.\\
     
       Nous avons clairement~:
      
           \begin{remarque}  Soit $M$  un ${\cal D}_X[s_1, \ldots ,s_p] $-Module coh\'erent, pour tout $c \in {\bf C}^p$~:
 $ {\rm car}^{\sharp} \,  M(c)   \subset  ({\rm car}^{\sharp} \, M)(0)$  et $ {\rm car}^{\rm rel} \,  M(c)   \subset ({\rm car}^{\rm rel} \, M)(c)$. Nous identifions ces sous-espaces \`a des sous-ensembles de  $T^{\ast} X$.
         
   \end{remarque}
       
        Soit localement, $m_1, \ldots ,m_l$ un syst\`eme de g\'en\'erateurs de $M$. Posons pour  tout entier $k$ :
       $$ M(c) (k) = \frac{\sum_{j=1}^l {\cal D}_X(k)m_j + \sum_{i=1}^p  (s_i -c_i )M}{ \sum_{i=1}^p  (s_i -c_i )M} \; .$$
       Les $M(c) (k) $ sont  \`a la fois une bonne filtration di\`ese,   une bonne filtration  relative de  $M(c) $, mais aussi une bonne filtration de $M(c)$ comme
       ${\cal D}_X$-Module. L'id\'eal engendr\'e par  $s_1,\ldots ,s_p$ (resp. $(s_1 - c_1,\ldots ,s_p - c_p)$)  et par  les symboles des op\'erateurs de ${\cal D}_X$ annulant  $ {\rm gr}\, ( M(c)) $ est l'id\'eal des symboles di\`eses (resp. relatifs) des op\'erateurs de ${\cal D}_X[s_1, \ldots ,s_p] $  annulant  $ {\rm gr}\, ( M(c) )$.
      Nous en d\'eduisons avec les identifications \'evidentes~:
       
       \begin{remarque}  Soit $M$  un ${\cal D}_X[s_1, \ldots ,s_p] $-Module coh\'erent, pour tout $c \in {\bf C}^p$~:
       $$ {\rm car}_{ {\cal D}_X  }\,   M(c)  = {\rm car}^{\sharp} \, M(c)   = {\rm car}^{\rm rel} \,  M(c)  \; . $$ 
   \end{remarque}
    En particulier, nous obtenons :
   \begin{proposition} Soit $M$   un ${\cal D}_X[s_1, \ldots ,s_p] $-Module coh\'erent major\'e  par une lagrangienne $ \Lambda $, 
   alors pour tout   $c \in {\bf C}^p$, les  ${\cal D}_X$-Modules $M(c)$ sont holonomes de vari\'et\'es caract\'eristiques contenues dans $ \Lambda $.
   \end{proposition}
       
       Nous allons pr\'eciser maintenant la structure de la vari\'et\'e caract\'eristique relative d'un ${\cal D}_X[s_1, \ldots ,s_p] $-Module coh\'erent major\'e  par une lagrangienne.

       \begin{proposition} \label{psvcrmml} Soit $M$  un ${\cal D}_X[s_1, \ldots ,s_p] $-Module coh\'erent   major\'e  par une 
  lagrangienne  de composantes irr\'eductibles   $( T_{X_{\alpha}}^{\ast} X )_{\alpha \in A}$, alors il existe des sous-vari\'et\'es alg\'ebriques $S_{\alpha}$ de
 ${\bf C}^p$ telles que~:
 $$  {\rm car}^{\rm rel} \, M = \bigcup \;   T_{X_{\alpha} }^{\ast} X  \times  S_{\alpha}   \; .$$  
       \end{proposition}
 
  \noindent {\bf Preuve :  }  Par hypoth\`ese, pour tout $ c \in {\bf C}^p$, $({\rm car}^{\rm rel} \, M)(c)$ est  contenu dans la r\'eunion des  
   $T_{X_{\alpha}}^{\ast} X $. D'apr\`es le th\'eor\`eme d'involutivit\'e, ses composantes irr\'eductibles sont de dimensions sup\'erieures \`a la dimension de $X$.
      Il en r\'esulte que pour tout  $ c \in {\bf C}^p$, les composantes irr\'eductibles de  $({\rm car}^{\rm rel} \, M)(c)$ sont certaines des lagrangiennes    $T_{X_{\alpha}}^{\ast} X $.  Il existe donc des sous-ensembles $S_{\alpha} $   de $ {\bf C}^p$ tels que :
    $$  {\rm car}^{\rm rel} \, M = \bigcup \,   T_{X_{\alpha} }^{\ast} X  \times  S_{\alpha}   \; .$$ 
    Il reste \`a montrer que les    $S_{\alpha} $ sont des sous-vari\'et\'es alg\'ebriques   de
 ${\bf C}^p$. Soit $( x_{\alpha}, \xi _{\alpha}) $  un point g\'en\'erique de $T_{ X_{\alpha} }^{\ast} X$, nous avons~: 
 $$ S_{\alpha}  =   ({\rm car}^{\rm rel} \, M ) \bigcap \, \{ (x , \xi, s)  \in  T^{\ast} X  \times  {\bf C}^p  \;  ;  \;  x = x_{\alpha}  \;  { \rm et }  \;  \xi = \xi _{\alpha}  \}  \; .$$
 Comme ${\rm car}^{\rm rel} \, M$ est localement d\'efini par un nombre fini d'\'equations polynomiales en $s$,   les  $S_{\alpha} $ sont des sous-vari\'et\'es alg\'ebriques   de  ${\bf C}^p$.

   
   \begin{definition} \label{didbg} Soit $M$  un ${\cal D}_X[s_1, \ldots ,s_p] $-Module coh\'erent  et $x_0 \in X$. Nous appelons id\'eal de Bernstein de
   $M$ en $x_0$ l'id\'eal  not\'e ${\cal B}_{x_0}(M) $ de $ {\bf C}[s_1, \ldots ,s_p] $ d\'efini par~:
   $$ {\cal B}_{x_0}(M) = \{ b(s_1, \ldots ,s_p   ) \in {\bf C}[s_1, \ldots ,s_p] \; ; \; b(s_1, \ldots ,s_p   ) M_{x_0} = 0  \}  \; .$$
   \end{definition}
   
  \begin{proposition}  \label{pvcrib} Soit $M$   un ${\cal D}_X[s_1, \ldots ,s_p] $-Module coh\'erent  de vari\'et\'e caract\'eristique relative  
    $\bigcup_{\alpha \in A } \,   T_{X_{\alpha} }^{\ast} X  \times  S_{\alpha} $ o\`u $ X_{\alpha}$ (resp . $S_{\alpha}$)
    sont des sous-espaces analytiques de $X$ (resp. sous-vari\'et\'es alg\'ebriques  de
 ${\bf C}^p$), Alors   $V({\cal B}_{x_0}(M)  )$  la vari\'et\'e des z\'eros de l'id\'eal de Bernstein de
   $M$ en $x_0$ est~:
   $$ V({\cal B}_{x_0}(M)  )  =  \bigcup_{\alpha \in A, x_0 \in  X_{\alpha} } \, S_{\alpha} \; .$$
   \end{proposition}
   
    \noindent {\bf Preuve :  }  Soit  $b(s_1, \ldots ,s_p   ) \in {\cal B}_{x_0}(M)$. Nous  pouvons voir $b$ comme un op\'erateur diff\'erentiel relatif de degr\'e z\'ero.
    Son symbole relatif  co\"{\i}ncide avec lui m\^eme. Il
     s'annule  d'apr\`es l'hypoth\`ese sur  ${\rm car}^{\rm rel} \, M$. Ainsi , $b(s_1, \ldots ,s_p   )$ s'annule sur les $S_{\alpha} $
     telles que $x_0 \in  X_{\alpha} $. Inversement, supposons que  $b(s_1, \ldots ,s_p   )$ s'annule sur les $S_{\alpha} $
     telles que $x_0 \in  X_{\alpha} $. Il s'annule   donc sur  ${\rm car}^{\rm rel} \, M$  au voisinage de $x_0$.  D'apr\`es le corollaire  \ref{ctzr} du  th\'eor\`eme des z\'eros relatif ,  il existe un entier $k$
     tel que  $b^k \in {\cal J}^{\rm rel}\, (M)$ au voisinage de $x_0$. Nous en d\'eduisons pour toute section $m$ de $M_{x_0} $, l'existence d'un entier $k$ tel que
     $b^k m  =0$ . Il en r\'esulte   que pour $k$ assez grand       $b^k \in {\cal B}_{x_0}(M) $ et  que $b$ s'annule sur $V({\cal B}_{x_0}(M)  ) $.
     
     \begin{corollaire}\label{cib} Soit $M$ un ${\cal D}_X[s_1, \ldots ,s_p] $-Module coh\'erent major\'e par une lagrangienne et $x_0 \in X$.
     Si ${\rm grade}\, M_{x_0}  \geq {\rm dim} \, X +1$, l'id\'eal de Bernstein de $M$ en $x_0$ est non r\'eduit \`a z\'ero.
     \end{corollaire}
   
       \noindent {\bf Preuve :  } Suivant la remarque \ref{ridg}, pour tout $a \in  {\bf C}^p$~:
       $${\rm dim }_{x_0,a,0} {\rm car}^{\rm rel} \, M \leq 2 \, {\rm dim } \, X +p -   {\rm grade}\, ( {\rm gr }^{\rm rel} \, M )_{x_0}  \leq
       {\rm dim } \, X +p -1 \; .$$
       Or suivant  la proposition \ref{psvcrmml},  il existe des sous-vari\'et\'es alg\'ebriques $S_{\alpha}$ de
 ${\bf C}^p$ telles que~:
 $$  {\rm car}^{\rm rel} \, M = \bigcup\,   T_{X_{\alpha} }^{\ast} X  \times  S_{\alpha}   \; .$$  
 Donc, pour $x_0 \in X_{\alpha} $,  la vari\'et\'e  alg\'ebrique  $S_{\alpha}$ est  donc de dimension strictement inf\'erieure \`a $p$.
Il reste  \`a utiliser la proposition \ref{pvcrib}.

\subsection{${\cal D}_X[s_1, \ldots ,s_p] $-Module major\'e par une lagrangienne sans ${\bf C}[s_1, \ldots ,s_p] $-torsion}

\begin{definition} 
Soit $B$ un anneau  r\'egulier qui v\'erifie la condition d'Auslander. Soit $M$ un $B$-module,  nous disons que 
$M$ est pur si tous les sous-modules de $M$ non r\'eduit \`a z\'ero ont m\^eme nombre grade.
\end{definition}

\begin{proposition}\label{pstep} Soit $M$ un ${\cal D}_X[s_1, \ldots ,s_p] $-Module coh\'erent major\'e  par une lagrangienne. Les assertions 
suivantes sont \'equivalentes~: 
\begin{enumerate}
\item Les sections non nulles des  fibres    de  $M$ sont sans  ${\bf C}[s_1, \ldots ,s_p] $-torsion.
\item  Les fibres non nulles de $M$   sont des modules   purs de nombre grade ${\rm dim}\, X$.
\end{enumerate}
De plus, si l'une de ces assertions est v\'erifi\'ee, il existe une vari\'et\'e lagrangienne $\Lambda$ de 
$T^{\ast} X$ telle que ${\rm car}^{\rm rel} \, M = \Lambda \times {\bf C}^p$.
    \end{proposition}

     \noindent {\bf Preuve}  $1 \Rightarrow 2$ :   D'apr\`es la proposition \ref{psvcrmml} ,  il existe des sous-espaces analytiques (resp. alg\'ebriques)
     $X_{\alpha}$ de $X$  (resp.  $S_{\alpha}$ de
 ${\bf C}^p$) telles que~:
 $$  {\rm car}^{\rm rel} \, M = \bigcup\,   T_{X_{\alpha} }^{\ast} X  \times  S_{\alpha}   \; .$$     
 Soit $x_0 \in X$, suivant la proposition \ref{pgfr}~:
 $$ {\rm grade} \,  ({\rm gr }^{\rm rel}\, M)_{x_0} = 2 \, {\rm dim } \, X +p - sup \{ {\rm dim} \, S_{\alpha} \; ; \;  x_0 \in  X_{\alpha} \} \; .$$
 Il en r\'esulte~: $ {\rm grade}\,  ({\rm gr}^{\rm rel}\, M)_{x_0} \geq {\rm dim } \, X$. Rappelons que suivant la proposition \ref{pegr},
 ${\rm grade  }  \, M_{x_0}     = {\rm grade  }\,    ( {\rm gr }^{\rm rel} \, M)_{x_0} $. Soit $N$ un sous-module de $M_{x_0}$, l'anneau 
 ${\cal D}_{X,x_0}[s_1, \ldots ,s_p] $ \'etant r\'egulier et v\'erifiant la condition d'Auslander : $ {\rm grade}\, N \geq {\rm grade  }  \, M_{x_0}   $
 (\cite{Bj2} , appendice IV proposition 2.3 ). Le module $N$ est major\'ee par une lagrangienne comme sous-module d'un module
    major\'ee par une lagrangienne.
 Si   ${\rm grade}\, N   \geq {\rm dim} \, X +1 $, il r\'esulterait du corollaire \ref{cib}
 que $N$ aurait de la ${\bf C}[s_1, \ldots ,s_p] $-torsion ce qui est contraire aux hypoth\`eses. Ainsi,  $ M_{x_0}$  est pur de  grade ${\rm dim}\, X$.\\
 
     \noindent {\bf Preuve}  $2 \Rightarrow 1$ : Si pour $x_0 \in X$, une section non nulle $M_{x_0}$ avait de la ${\bf C}[s_1, \ldots ,s_p] $-torsion, $M$ admettrait au voisinage de $x_0$
    un  sous-module $N$ tel que  ${\cal B}_{x_0}(N)\neq 0$.  Il  r\'esulterait des propositions 
 \ref{psvcrmml} et  \ref{pvcrib} que  les vari\'et\'es alg\'ebriques $  S_{\alpha}$ pour  $x_0 \in  X_{\alpha}  $ intervenant dans la vari\'et\'e caract\'eristique
 relative de $N$ seraient  toutes de dimensions inf\'erieures ou \'egales \`a $p-1$. Il r\'esulterait de la proposition \ref{pgfr}  que
 ${\rm grade}\, N   \geq {\rm dim} \, X +1 $ ce qui contredit l'hypoth\`ese de puret\'e.\\
 
    \noindent {\bf Fin de la preuve}    :  Supposons v\'erifi\'ees ces conditions \'equivalentes. Suivant \cite{Bj2} Appendice IV  theorem 4.11,  il existe une bonne filtration relative  $\Gamma$ de $M_{x_0}$ dont le gradu\'e $ {\rm gr}_{\Gamma}\, M_{x_0}$ est un ${\cal D}_{X,x_0}[s_1, \ldots ,s_p] $ pur de  grade ${\rm dim}\, X$. Rappelons,  \cite{Bj2}, Appendice IV,  proposition 2.6, que si $B$ est un anneau r\'egulier v\'erifiant la condition d'Auslander, un $B$ module $E$ de type fini est pur  de  grade $l$ si et seulement si $Ext_A^j(Ext_A^j(E,A),A) = 0$ \'equivaut \`a  $j \neq l$. Il en r\'esulte que pour $j\neq {\rm dim}\, X$~:
 $$ Ext ^j(Ext ^j({\rm gr} M, {\rm gr}\,{\cal D}_{X,}[s_1, \ldots ,s_p]),{\rm gr}\,{\cal D}_{X}[s_1, \ldots ,s_p])_{x_0}= 0 \; .$$ 
 En analytisant pour la filtration d\'eduite de $\Gamma$, nous obtenons pour tout $a \in { \bf C}^p$ et $j\neq {\rm dim}\, X$~:
 $$  Ext_{{\rm gr}\, {\cal D}_{X \times { \bf C}^p / { \bf C}^p  } }^j(Ext_{{\rm gr}\, {\cal D}_{X \times { \bf C}^p / { \bf C}^p   } }^j(  {\rm gr}\, M^{\rm an}, {\rm gr}\, {\cal D}_{X \times { \bf C}^p / { \bf C}^p    }),{\rm gr}\, {\cal D}_{X \times { \bf C}^p / { \bf C}^p    })_{x_0,a} = 0 \; .$$  
 Ainsi, les  fibres    $ ({\rm gr} \, M^{\rm an})_{x_0,a}$ non nulles  sont pures de  grade ${\rm dim}\, X$. 
 Rappelons que sur un anneau commutatif r\'egulier, si un   module de type fini  est pur de grade $l$,  tous les id\'eaux premiers minimaux de son support   sont  de  hauteur $l$
 (\cite{Bj2} Appendice IV  proposition 3.7). Nous obtenons ansi que les composantes irr\'eductibles de
 $ {\rm car}^{\rm rel} \, M $ au voisinage de $x_0$ sont de dimension ${\rm dim}\, X + p$. Et pour tout $\alpha$ tel que $x_0 \in  X_{\alpha}  $,  
 $S_{\alpha} = { \bf C}^p$.\\

     \begin{proposition} \label{pcrf} Soit $M$ un ${\cal D}_X[s_1, \ldots ,s_p] $-Module coh\'erent tel que ${\rm car}^{\rm rel} \, M = \Lambda \times {\bf C}^p$
     o\`u $\Lambda$ est une vari\'et\'e lagrangienne de $T^{\ast} X $. Alors, pour tout $c \in { \bf C}^p$~:
     $$ {\rm car}_{ {\cal D}_X  }\, (M(c))  = \Lambda \; .$$
    \end{proposition}

    \noindent {\bf Preuve} :  Il r\'esulte du fait que $   {\rm car}^{\rm rel} \, M  = \Lambda \times  {\bf C}^p $ que ${\cal J}^{\rm rel} (M) = I[s_1, \ldots ,s_p]   $ l'id\'eal de ${\cal O}_X[\xi _1, \ldots ,\xi _n, s_1, \ldots ,s_p]$ 
engendr\'e par l'id\'eal   
 $I$   de
${\cal O}_X[\xi _1, \ldots ,\xi _n]$ des fonctions nulles sur $\Lambda $. L'inclusion ${\cal J}^{\rm rel} (M) \subset I[s_1, \ldots ,s_p]  $ est claire. L'autre inclusion r\'esulte du
th\'eo\`eme des z\'eros relatifs,   proposition \ref{thzr} modulo  l'astuce de Rabinowitch. Soit $p_i$ pour $i = 1,\ldots ,l$, les id\'eaux premiers minimaux contenant
$I_{x_0}$, les $p_i[s_1, \ldots ,s_p]$ sont alors les id\'eaux premiers minimaux de  ${\cal J}^{\rm rel} (M)_{x_0} $.  
Consid\'erons le complexe \`a deux termes de multiplication par $s_1 -c_1$~:
$$  {\rm gr }^{\rm rel} \, M_{x_0} \;   \stackrel{s_1 - c_1}{\longrightarrow}   \; {\rm gr }^{\rm rel} \, M_{x_0} \; .$$
Notons   $L_1$  (resp. $K_1$) son conoyau (resp. noyau). Soit $ ({\cal J}^{\rm rel} (M)_{x_0}, s_1 -c_1)$ l'id\'eal de $
{\rm gr }^{\rm rel} \, {\cal D}_{X,x_0}[s_1, \ldots ,s_p] $ engendr\'e par   ${\cal J}^{\rm rel} (M)_{x_0}$ et $ s_1 -c_1$.
Ses id\'eaux premiers minimaux sont les $q_i = (  p_i[s_1, \ldots ,s_p] , s_1 -c_1 ) $ pour $i = 1,\ldots ,l$. Les localis\'es
$(L_1)_{q_i}$ et $(K_1)_{q_i}$ sont de longueur finie. Notons pour tout id\'eal premier $q$ d'un anneau commutatif $A$ ,  $e_{q}(E)$ la multiplicit\'e d'un
$A$ module $E$ de type fini. Un calcul de multiplicit\'e donne~:
$$ e_{q_i}(L_1) - e_{q_i }(K_1)  = e_{p_i}( {\rm gr}^{\rm rel} \, M) $$
De plus, M \'etant muni d'une bonne filtration le complexe :
$$M \; \stackrel{s_1 - c_1}{\longrightarrow}   \; M $$ de multiplication par
$s_1 - c_1$ est filtr\'e. Le fait que cette filtration soit bonne implique la convergence de la suite spectrale de ce complexe filtr\'e.
La filtration de $M$ induit sur $L_1$ et $K_1$ des filtrations \`a gradu\'es noeth\'eriens, donc de bonnes filtrations. Nous obtenons~:
$$ e_{q_i }( {\rm gr}^{\rm rel} \, L_1 ) - e_{q_i }( {\rm gr}^{\rm rel} \, K_1 )  = e_{p_i }( {\rm gr}^{\rm rel} \, M ) > 0 $$
Il en r\'esulte~:
$$ e_{q_i }( {\rm gr}^{\rm rel} \, L_1) > 0\; .$$
Ainsi, les $q_i$ sont des id\'eaux premiers minimaux du support de $ {\rm gr}^{\rm rel} \, L_1$. D'autre part,
$( {\cal J}^{\rm rel} (M)_{x_0} , s_1 - c_1 ) \subset  {\cal J}^{\rm rel} (L_1)$. Donc, chaque  $q_i$ est contenu dans un id\'eal
premier minimal du support de   
$ {\rm gr}^{\rm rel} \, L_1$.
 Il en r\'esulte que
 les $q_i$ sont les  id\'eaux premiers minimaux
de ${\rm gr}^{\rm rel} \, L_1 $. On obtient ainsi,
$${\rm car}^{\rm rel} \, L_1 = \Lambda \times \{c_1\}\times {\bf C}^{p-1} \; .$$  Il reste
\`a it\'erer pour obtenir $  {\rm car}^{\rm rel} \,( M(c)) = \Lambda$.

\begin{proposition}\label{pcvcf}  Soit $M$ un ${\cal D}_X[s_1, \ldots ,s_p] $-Module coh\'erent dont les fibres sont des modules  purs de grade ${\rm dim } \, X $ et   tel que $({\rm car}^{\sharp} \, M )(0)$ soit lagrangienne.
 Alors, pour tout
$c \in {\bf C}^p $ :
$$({\rm car}^{\sharp} \, M )(0) =  {\rm car}^{\sharp} \,(  M(c) )  =
{\rm car}_{ {\cal D}_X }   \,( M )(c)   $$
\end{proposition}

\noindent {\bf Preuve : } Il faut montrer que $({\rm car}^{\sharp} \, M )(0) =  {\rm car}^{\sharp} \,(  M(0) )$.
Par hypoth\`ese, $({\rm car}^{\sharp} \, M )(0)$ est lagrangienne et toutes les composantes irr\'eductibles de ${\rm car}^{\sharp} \, M $
sont de dimension ${\rm dim } \, X +p$. Il en r\'esulte qu'aucune composante de ${\rm car}^{\sharp} \, M $ n'est contenue dans $s_i=0$. Filtrons $M$ par une bonne filtration di\`ese.
La proposition s'obtient par r\'ecurrence en consid\'erant la suite spectrale du complexe filtr\'e de multiplication par $s_1 - c_1$~:
$$M \; \stackrel{s_1 - c_1}{\longrightarrow}   \; M \; $$


   \newpage

   \section{Equations fonctionnelles associ\'ees \`a Module holonome relativement \`a plusieurs fonctions analytiques}
   \label{s3}
 
Soit $f_1, \ldots, f_p$, des fonctions analytiques sur $X$. Notons $F= f_1f_2 \cdots f_p$ leur produit.
   Soit $M$ un ${ \cal D}_X$-Module holonome. Consid\'erons ${\cal O}_X[s_1, \ldots ,s_p, 1/F] f_1^{s_1} \ldots  f_p^{s_p} $  le ${\cal O}_ X [s_1, \ldots ,s_p, 1/F]$-Module libre de rang $1$ de base $f_1^{s_1} \ldots  f_p^{s_p} $. Le produit tensoriel :
$$ M \otimes_{ {\cal O}_ X} {\cal O}_X [s_1, \ldots ,s_p,1/ F] f_1^{s_1} \ldots  f_p^{s_p}$$ est muni d'une structure naturelle
de ${\cal D}_X [s_1, \ldots ,s_p]$-Module d\'efinie pour toute  section  $ n$ de $M$ et $  a $  de 
${\cal O}_X[s_1, \ldots ,s_p,1/ F]$ par~:
$$ \frac{\partial }{\partial x_i } (n \otimes a  f_1^{s_1} \ldots  f_p^{s_p} ) = $$
$$
 ( \frac{\partial }{\partial x_i } n  ) \otimes a  f_1^{s_1} \ldots  f_p^{s_p} + 
n  \otimes \frac{\partial a}{\partial x_i }f_1^{s_1} \ldots  f_p^{s_p} + \sum_{j=1}^p
 n   \otimes s_j a   \frac{ \frac{\partial f_j}{\partial x_i }}{f_j} f_1^{s_1} \ldots  f_p^{s_p} \quad   .$$
 Nous notons, $   n  \otimes a  f_1^{s_1} \ldots  f_p^{s_p} = a\, n f_1^{s_1} \ldots  f_p^{s_p}  $ et
 $f_1^{s_1+r} \ldots  f_p^{s_p+r} = F^r f_1^{s_1} \ldots  f_p^{s_p}$ pour tout entier $r$.\\
 
 Nous supposons $M$ engendr\'e par une section $m$. \\

   \subsection{Vari\'et\'e carat\'eristique relative de    ${\cal D}_X[s_1, \ldots ,s_p] m f_1^{s_1} \ldots  f_p^{s_p}$}

  Dans ce paragraphe, $M$ d\'esigne un module holonome engendr\'e par une section  $m$. En dehors de $F=0$, pour toute section $m$ de $M$~: 
 $$ \frac{\partial }{\partial x_i } (n \otimes    f_1^{s_1} \ldots  f_p^{s_p} ) =  
\left[   ( \frac{\partial }{\partial x_i } + \sum_{i=1}^p \frac{s_i}{f_i}  \frac{\partial f_i}{\partial x_i }  ) n        \right]
 f_1^{s_1} \ldots  f_p^{s_p}  .$$
 Munissons $M$ de la bonne filtration $ {\cal D}_X(k)m$ et  ${\cal D}_X[s_1, \ldots ,s_p] m f_1^{s_1} \ldots  f_p^{s_p} $
 de la bonne filtration relative 
 ${\cal D}_X[s_1, \ldots ,s_p](k) m f_1^{s_1} \ldots  f_p^{s_p} $. En dehors de $F=0$, nous obtenons~:
 $$ {\rm gr }^{\rm rel}\,  {\cal D}_X[s_1, \ldots ,s_p] m f_1^{s_1} \ldots  f_p^{s_p}  \simeq
( {\rm gr }_{{\cal D}_X}\, M) [s_1, \ldots ,s_p] \; .$$
Comme $M$ est holonome, pour tout $x_0$ dans le support de $M$, le grade de $({\rm gr }_{{\cal D}_X}\, M)_{x_0}$
est \'egal \`a ${\rm dim}\, X$. Il en r\'esulte  que si $F(x_0)\neq 0$~: 
$$        {\rm grade}\, ( {\rm gr }^{\rm rel}\,  {\cal D}_X[s_1, \ldots ,s_p] m f_1^{s_1} \ldots  f_p^{s_p} )_{x_0}
= {\rm grade}\, ( {\rm gr }_{{\cal D}_X}\, M)_{x_0}  = {\rm dim}\, X \; .$$
Nous obtenons ainsi~:
\begin{lemme} Pour tout $x_0$ dans le support de $M$ tel que $F(x_0)\neq 0$~:
$$ {\rm grade}\, ( {\rm gr }^{\rm rel}\,  {\cal D}_X[s_1, \ldots ,s_p] m f_1^{s_1} \ldots  f_p^{s_p} )_{x_0}
 =   {\rm dim}\, X \; .$$
\end{lemme}

 Notons ${\cal G}_k( {\cal D}_X[s_1, \ldots ,s_p] m f_1^{s_1} \ldots  f_p^{s_p} ) $ le sous-Module form\'e des sections  
 de ${\cal D}_X[s_1, \ldots ,s_p] m f_1^{s_1} \ldots  f_p^{s_p}$   dont les fibres sont de grade sup\'erieure ou \'egal \`a $k$. Les faisceaux
 $ {\rm gr}^{\rm rel}\,  {\cal D}_X[s_1, \ldots ,s_p]$ et ${\cal D}_X[s_1, \ldots ,s_p]$ v\'erifient la condition de Noether~:  sur tout espace localement compact, les suites croissantes de sous-Modules coh\'erents stationnent. Cela permet d'assurer (\cite{Bj2} Appendice IV theorem 2.30) 
  que pour tout entier $k$, les  ${\cal G}_k( {\cal D}_X[s_1, \ldots ,s_p] m f_1^{s_1} \ldots  f_p^{s_p} ) $ sont des 
  ${\cal D}_X[s_1, \ldots ,s_p]$-modules coh\'erents. En effet, le gradu\'e de cette filtration ${\cal G}_.$ s'obtient comme limite d'une suite spectrale
  associ\'ee \`a un complexe double filtr\'e dont le deuxi\`eme tableau est compos\'e des ${\cal D}_X[s_1, \ldots ,s_p]$-modules coh\'erents~:
  $$Ext_{{\cal D}_X[s]}^j(Ext_{{\cal D}_X[s]}^j( {\cal D}_X[s_1, \ldots ,s_p] m f_1^{s_1} \ldots  f_p^{s_p} , {\cal D}_X[s_1, \ldots ,s_p]),  {\cal D}_X[s_1, \ldots ,s_p])   \, . $$
  
\begin{lemme} \label{lmgef} Pour tout $x_0$ dans le support de $M[1/F]$~:
$$ {\rm grade}\, ( {\rm gr }^{\rm rel}\,  {\cal D}_X[s_1, \ldots ,s_p] m f_1^{s_1} \ldots  f_p^{s_p} )_{x_0}
\geq {\rm dim}\, X \; .$$
\end{lemme}

\noindent {\bf Preuve : } Le quotient~:
$$\frac{{\cal D}_X[s_1, \ldots ,s_p] m f_1^{s_1} \ldots  f_p^{s_p}}{{\cal G}_{ {\rm dim}\, X}( {\cal D}_X[s_1, \ldots ,s_p] m f_1^{s_1} \ldots  f_p^{s_p}} \; .
$$  
est un ${\cal D}_X[s_1, \ldots ,s_p]$ Module coh\'erent support\'e par $F=0$. Du th\'eor\`eme des z\'eros,     au voisinage de tout $x_0\in X$,
il r\'esulte l'existence d'un entier $l$ tel que~:
$$ m F^lf_1^{s_1} \ldots  f_p^{s_p} = m  f_1^{s_1+ l} \ldots  f_p^{s_p + l}  \in {\cal G}_{ {\rm dim}\, X}( {\cal D}_X[s_1, \ldots ,s_p] m f_1^{s_1} \ldots  f_p^{s_p})\; .
$$  
Comme pour tout $P (x, \partial _x,s) \in {\cal D}_X[s_1, \ldots ,s_p]$~:
$$ P (x, \partial _x,s)  m f_1^{s_1} \ldots  f_p^{s_p} = 0 \Longleftrightarrow  P (x, \partial _x,s_1+l, \ldots,s_p+l)  mf_1^{s_1+l} \ldots  f_p^{s_p+l} =0\, ,$$
nous d\'eduisons d'une r\'esolution libre locale de  ${\cal D}_X[s] m f_1^{s_1} \ldots  f_p^{s_p}
$ une r\'esolution libre locale de ${\cal D}_X[s] m f_1^{s_1+l} \ldots  f_p^{s_p+l}$ et les bijections de translations entre~:
$$Ext^i({\cal D}_X[s] m f_1^{s_1} \ldots  f_p^{s_p}, {\cal D}_X[s] )\; {\rm et } \;   
Ext^i({\cal D}_X[s] m f_1^{s_1+l} \ldots  f_p^{s_p+l}, {\cal D}_X[s] )  \; .$$
Il en r\'esulte que les  fibres en $x_0$ de ${\cal D}_X[s] m f_1^{s_1} \ldots  f_p^{s_p}$ et 
 ${\cal D}_X[s] m f_1^{s_1+l} \ldots  f_p^{s_p+l}$ ont m\^eme  grade et donc le lemme.
 
\begin{proposition}  \label{pefml} II existe une vari\'et\'e lagrangienne conique  $\Lambda$ de $T^{\ast}X$ tel que~:
$$ {\rm car}^{\rm rel} \,  {\cal D}_X[s] m f_1^{s_1} \ldots  f_p^{s_p} = \Lambda \times {\bf C}^p \; .$$
\end{proposition}

 \noindent {\bf Preuve : }  Commen\c{c}ons par donner une preuve g\'eom\'etrique.
 Noton $1_i= (0\ldots, 0,1,0\ldots0)$ o\`u le $1$ est pla\c{c}\'e \`a  la i-\`eme place et appelons $\tau _i$
 la translation  de $T^{\ast}X \times {\bf C}^p$ d\'efinie par    $(x,\xi,s) \mapsto  (x,\xi,s+1_i)$. Nous avons 
 pour tout $P (x, \partial _x,s) \in {\cal D}_X[s_1, \ldots ,s_p]$~:
$$ P (x, \partial _x,s)  m f_1^{s_1} \ldots  f_p^{s_p} = 0 \Longleftrightarrow  P (x, \partial _x,s+ 1_i)  mf_if_1^{s_1} \ldots  f_p^{s_p} =0\, .$$
Ainsi,
$${\rm car}^{\rm rel} \, (    {\cal D}_X[s_1, \ldots ,s_p] mf_if_1^{s_1} \ldots  f_p^{s_p}  ) = {\tau _i}^{-1} (
{\rm car}^{\rm rel} \, (    {\cal D}_X[s_1, \ldots ,s_p] m f_1^{s_1} \ldots  f_p^{s_p}  )  \; .$$

 D'autre part,  il est clair que $ {\cal D}_X[s_1, \ldots ,s_p] mf_if_1^{s_1} \ldots  f_p^{s_p} $ est un sous-Module de $ {\cal D}_X[s_1, \ldots ,s_p] mf_1^{s_1} \ldots  f_p^{s_p} $. Donc, ${\rm car}^{\rm rel} \, (    {\cal D}_X[s_1, \ldots ,s_p] m f_1^{s_1} \ldots  f_p^{s_p}  )  $ est stable par les translations
 $  {\tau _i}^{-1}$ et donc par toutes translations $(x,\xi,s)  \mapsto (x,\xi,s-k)$ de $T^{\ast}X \times {\bf C}^p$ o\`u $k \in {\bf N}^p$. Nous en d\'eduisons 
l'existence  d'un sous-ensemble analytique conique $\Lambda$ de $T^{\ast}X $ tel que~:
$$ {\rm car}^{\rm rel} \,  {\cal D}_X[s] m f_1^{s_1} \ldots  f_p^{s_p} = \Lambda \times {\bf C}^p \; .$$
D'apr\`es le lemme \ref{lmgef} et la remarque \ref{ridg} de la section \ref{ssfr}, pour tout $x_0$ du support de $M[1/F]$ et $a\in {\bf C}^p$~:
$$ {\rm dim }_{x_0,a,0} \, {\rm car}^{\rm rel} \,  {\cal D}_X[s] m f_1^{s_1} \ldots  f_p^{s_p} \leq  {\rm dim }  \, X+ p\; .$$
Donc, ${\rm dim }_{x_0,0} \, \Lambda \leq n$. Comme $\Lambda$ est conique,  la dimension de  $\Lambda$ en tout point est inf\'erieure
\`a $n$.  En tout  point lisse de $\Lambda$, la projection de $ \Lambda \times {\bf C}^p $ sur  $  {\bf C}^p$ est une submersion. Il r\'esulte du th\'eor\`eme d'involutivit\'e
  rappel\'e au paragraphe \ref{scti}   que $\Lambda$ est un sous espace involutif conique de $T^{\ast}X $. C'est donc une vari\'et\'e lagrangienne.\\

Donnons une preuve plus alg\'ebrique. Soit $P(s) \in {\cal D}_X^{\rm rel }[s_1, \ldots ,s_p](k)$ tel que $P(s)m f_1^{s_1} \ldots  f_p^{s_p}=0$.
Remarquons que pour tout $r\in {\bf N}^p$:
$$\begin{array}{rcl} P(s-r)m f_1^{s_1} \ldots  f_p^{s_p} &=& P(s-r)m f^r f_1^{s_1-r_1} \ldots  f_p^{s_p-r_p} \\
&=&  (P(s-r)f^r-f^rP(s-r)) m f_1^{s_1-r_1} \ldots  f_p^{s_p-r_p} \\
& \in & {\cal D}_X^{\rm rel }[s_1, \ldots ,s_p](k-1)m f_1^{s_1-r_1} \ldots  f_p^{s_p-r_p} \; .
\end{array}
 $$
 Ainsi, $\sigma^{\rm rel} (P(s-r))$ annule le gradu\'e de $ {\cal D}_X[s] m f_1^{s_1} \ldots  f_p^{s_p} $ pour la bonne filtration induite par la   filtration   naturelle de  ${\cal D}_X  [s_1, \ldots ,s_p] m f_1^{s_1-r_1} \ldots  f_p^{s_p-r_p}$. Donc,  $\sigma^{\rm rel} (P(s-r)) \in
 {\cal J}^{\rm rel }(  {\cal D}_X[s] m f_1^{s_1} \ldots  f_p^{s_p})$ pour tout  $ r\in {\bf N}^p$. Cela montre  que cet id\'eal est engendr\'e par des symboles
 ind\'ependants de $s$. Suivant le    th\'eor\'eme d'involutivit\'e sa racine est stable par crochet de Poisson.
 Nous en d\'eduisons~:
 $$ {\rm car}^{\rm rel} \,  {\cal D}_X[s] m f_1^{s_1} \ldots  f_p^{s_p} = \Lambda \times {\bf C}^p \; ,$$
o\`u $ \Lambda$ est une vari\'et\'e involutive conique de $T^{\ast}X $. Le m\^eme argument sur la dimension que pr\'ec\'edemment
utilisant le lemme \ref{lmgef} et la remarque \ref{ridg} de la section \ref{ssfr}   permet de conclure  que  $\Lambda$ est   lagrangienne.


\begin{proposition} \label{ppdmfs}
Pour tout $x_0$ dans le support de $M[1/F]$, le module $ {\cal D}_{X,x_0}[s_1, \ldots ,s_p] m f_1^{s_1} \ldots  f_p^{s_p}  $ est pur de grade 
${\rm dim}\, X$ : tout sous module non nul de $ {\cal D}_{X,x_0}[s_1, \ldots ,s_p] m f_1^{s_1} \ldots  f_p^{s_p} $  est de grade ${\rm dim}\, X$. 
 \end{proposition}

 \noindent {\bf Preuve : }  La proposition \ref{pefml} dit que le Module  $ {\cal D}_{X }[s_1, \ldots ,s_p] m f_1^{s_1} \ldots  f_p^{s_p}   $ est major\'ee par une lagrangienne. Or, $({\cal D}_X[s_1, \ldots ,s_p] m f_1^{s_1} \ldots  f_p^{s_p} )_{x_0}$ est sans 
 ${\bf C}[s_1, \ldots ,s_p]$-torsion. Le r\'esultat  r\'esulte de la proposition \ref{pstep}.

\begin{lemme} \label{lmgeff}  
Pour tout $x_0$ dans le support de $M[1/F]$  tel que $F(x_0)=0$ :
$$ {\rm grade}\,   \left( \frac {     {\cal D}_X[s_1, \ldots ,s_p] m f_1^{s_1} \ldots  f_p^{s_p}   }{   {\cal D}_X[s_1, \ldots ,s_p] m f_1^{s_1+1} \ldots  f_p^{s_p+1}  } \right)_{x_0}
\geq {\rm dim}\, X + 1\; .$$
\end{lemme}
 
 \noindent {\bf Preuve : } Soit $\Lambda$ 
la  vari\'et\'e lagrangienne conique   de $T^{\ast}X$ telle  que 
$$ {\rm car}^{\rm rel} \,  {\cal D}_X[s] m f_1^{s_1} \ldots  f_p^{s_p} = \Lambda \times {\bf C}^p \; .$$
L'id\'eal ${\cal J}^{\rm rel }(  {\cal D}_X[s] m f_1^{s_1} \ldots  f_p^{s_p})$ est l'id\'eal $I \otimes_ {{\rm gr}^{\rm rel }\,    {\cal D}_X} {\rm gr}^{\rm rel }\,    {\cal D}_X[s]$
o\`u $I$ est l'id\'eal de  ${\rm gr}^{\rm rel }\,    {\cal D}_X$ des fonctions nulles sur $\Lambda$. Soit $P$ un id\'eal premier d\'efinissant une composante irr\'eductible de 
$\Lambda$. Cet id\'eal $P$ est de dimension  $ {\rm dim}\, X$. Par un argument de translation, les multiplit\'es 
de ${\rm gr}^{\rm rel }\,  (  {\cal D}_X[s] m f_1^{s_1} \ldots  f_p^{s_p})_{x_0} $
 et  ${\rm gr}^{\rm rel }\,  (  {\cal D}_X[s] m f_1^{s_1+ 1} \ldots  f_p^{s_p+ 1})_{x_0} $ en   ${\cal P} \otimes_ {{\rm gr}^{\rm rel }\,    {\cal D}_X} {\rm gr}^{\rm rel }\,    {\cal D}_X[s]$ co\"{\i}ncident.
 Il en r\'esulte qu'aucune composant irr\'eductible de $\Lambda \times {\bf C}^p$ n 'est une composante de la vari\'et\'e caract\'eristique relative du quotient :
 $$\frac {     {\cal D}_X[s_1, \ldots ,s_p] m f_1^{s_1} \ldots  f_p^{s_p}   }{   {\cal D}_X[s_1, \ldots ,s_p] m f_1^{s_1+1} \ldots  f_p^{s_p+1}  }$$
 Ce quotient a donc une  vari\'et\'e  caract\'eristique  relative de dimension inf\'erieure ou \'egale \`a  $ {\rm dim}\, X +p - 1$. Il reste \`a utiliser la proposition \ref{pgfr}.


 \begin{definition} \label{dfnibm} Soit $m$ une section d'un ${ \cal D}_X$-Module holonome et $x_0 \in X$. Nous appelons id\'eal de Bernstein de 
  $(m, f_1, \ldots ,f_p) $ au voisinage de $x_0$ l'id\'eal ${\cal B}(m,x_0, f_1, \ldots ,f_p)$  de   ${\bf C}[s_1, \ldots ,s_p]$  form\'e des polyn\^omes   $b$ 
 v\'erifiant  au voisinage de $x_0$:
$$    b (s_1, \ldots ,s_p) m f_1^{s_1} \ldots  f_p^{s_p} \in {\cal D}_X[s_1, \ldots ,s_p] m f_1^{s_1+1} \ldots  f_p^{s_p+1} \; ,$$
qui sont appel\'es  polyn\^omes de Bernstein de $(m, f_1, \ldots ,f_p) $ au voisinage de $x_0$.
 
 \end{definition}

 Notons que l'id\'eal ${\cal B}(m,x_0, f_1, \ldots ,f_p)$ n'est autre que l'id\'eal de Bernstein, au sens de la d\'efinition \ref{didbg},  du 
${\cal D}_X[s_1, \ldots ,s_p] $-Module coh\'erent :
 $$  \frac{{\cal D}_X[s_1, \ldots ,s_p] m f_1^{s_1} \ldots  f_p^{s_p}}{{\cal D}_X[s_1, \ldots ,s_p] m f_1^{s_1+1} \ldots  f_p^{s_p+1}} \; .$$

  \begin{corollaire}   (Une preuve de l'existence d'une \'equation fonctionnelle de Bernstein dans le cas analytique) 
 Si $m$ est une  section d'un ${\cal D}_X $-Module holonome, au voisinage de tout $x_0 \in X$, il existe un polyn\^ome 
 $b(s_1, \ldots ,s_p)$ non nul~:
 $$  b(s_1, \ldots ,s_p) m f_1^{s_1} \ldots  f_p^{s_p} \in {\cal D}_X[s_1, \ldots ,s_p] m f_1^{s_1+1} \ldots  f_p^{s_p++1}\; . $$
 \end{corollaire}

  \noindent {\bf Preuve : } Par le lemme \ref{lmgef},  le quotient~:
   $$  \frac {     {\cal D}_X[s_1, \ldots ,s_p] m f_1^{s_1} \ldots  f_p^{s_p}   }{   {\cal D}_X[s_1, \ldots ,s_p] m f_1^{s_1+1} \ldots  f_p^{s_p+1}  }   $$
a ses fibres de de grade sup\'erieur ou \'egal \`a $ {\rm dim}\, X + 1$.   Il est major\'e par une lagrangienne, car quotient d'un  ${\cal D}_X[s_1, \ldots ,s_p]$-Module 
major\'e par une lagrangienne (proposition \ref{pefml}). Il reste \`a utiliser le corollaire \ref{cib}.

  \subsection{Construction d'un quotient pur d'un sous-facteur de ${\cal D}_X[s] m f^s$}
  
  Soit $m$ une section d'un ${ \cal D}_X$-Module holonome $M$.
  Consid\'erons les applications bijectives $\tau$~: 
  $$ M[\frac{1}{F} s_1, \ldots ,s_p]  f_1^{s_1} \ldots  f_p^{s_p} \stackrel{\tau}{\longrightarrow}
  M[\frac{1}{F} s_1, \ldots ,s_p]  f_1^{s_1} \ldots  f_p^{s_p}  \;\; {\rm et} \;\; {\cal D}_X  [s_1, \ldots ,s_p]   \stackrel{\tau}{\longrightarrow} {\cal D}_X  [s_1, \ldots ,s_p]   $$ 
  d\'efinies respectivement par~: 
$$ \tau ( n(s_1, \ldots ,s_p)f_1^{s_1} \ldots  f_p^{s_p} )=
  n(s_1+1, \ldots ,s_p+1)f_1^{s_1+1} \ldots  f_p^{s_p+1} $$
  et 
$$ \tau ( P(x,\partial _x, s_1, \ldots ,s_p) ) = P(x,\partial _x, s_1+1, \ldots ,s_p+p)   \; .$$ 
Nous d\'esignons \'egalement par $\tau$ l'application~:
$$ T^{\ast} X \times  {\bf C}^p  \stackrel{\tau}{\longrightarrow}  T^{\ast} X \times  {\bf C}^p \; \;
(x,\xi, s_1, \ldots ,s_p) \mapsto (x,\xi, s_1+1, \ldots ,s_p+1)\; .$$
 Pour toutes sections  $n$ de $M$ et tout $P \in {\cal D}_X  [s_1, \ldots ,s_p] $~:
 $$P \, n f_1^{s_1} \ldots  f_p^{s_p} = 0 \Longleftrightarrow \tau (P) \, n f_1^{s_1+1} \ldots  f_p^{s_p+1} = 0\; .$$
 Nous en d\'eduisons~:
 $$ {\rm car}^{\sharp} \,     {\cal D}_X[s_1, \ldots ,s_p] m f_1^{s_1+1} \ldots  f_p^{s_p+1}     = 
  {\rm car}^{\sharp} \,     {\cal D}_X[s_1, \ldots ,s_p] m f_1^{s_1} \ldots  f_p^{s_p}   \; ,$$
  $${\rm car}^{\rm rel} \,     {\cal D}_X[s_1, \ldots ,s_p] m f_1^{s_1+1} \ldots  f_p^{s_p+1}      =\tau ^{-1} ( {\rm car}^{\rm rel} \,  
 {\cal D}_X[s_1, \ldots ,s_p] m f_1^{s_1} \ldots  f_p^{s_p}   ) \; .$$

  L'objet de ce paragraphe est de d\'emontrer la proposition suivante et d'en donner les premi\`eres cons\'equences..

 \begin{proposition} \label{psfp} Soit $M$ un   holonome    engendr\'e par une section $m$. 
 Consid\'erons la
 famille ${\cal G}$  des     sous $ {\cal D}_{X}[s_1, \ldots ,s_p]$-Modules  $L$ de  type fini de $M[\frac{1}{F},s_1, \ldots ,s_p] f_1^{s_1 } \ldots  f_p^{s_p } $   contenant  $ {\cal D}_{X}[s_1, \ldots ,s_p] m f_1^{s_1 } \ldots  f_p^{s_p } $
  et tels que pour tout  tout point $x_0 \in X$ :
$$  {\rm grade} \, \frac{L_{x_0}}{ {\cal D}_{X,x_0}[s_1, \ldots ,s_p] m f_1^{s_1 } \ldots  f_p^{s_p } } \geq {\rm dim}\, X   +  2 \; $$
 Cette famille admet un plus grand \'el\'ement  not\'e   $\tilde{L}$. Ce Module $\tilde{L}$ est   un sous 
 ${\cal D}_{X}[s_1, \ldots ,s_p]$-Module coh\'erent de $M[\frac{1}{F},s_1, \ldots ,s_p] f_1^{s_1 } \ldots  f_p^{s_p } $
 v\'erifiant~:
\begin{enumerate}
\item ${\cal D}_{X}[s_1, \ldots ,s_p] m  f_1^{s_1} \ldots  f_p^{s_p} \subset \tilde{L}$,
\item $\tau ( \tilde{L} ) \subset \tilde{L} $,
\item $ \tilde{L} / \tau ( \tilde{L} )$ est un ${\cal D}_{X}[s_1, \ldots ,s_p]$-Module dont les fibres non nulles sont des modules purs de grade
${\rm dim}\, X +1$.
\end{enumerate}
\end{proposition}
 
  La preuve est une adaptation de celle du th\'eor\`eme 2.12 \cite{Bj2} appendice IV. Cette adaptation 
  est n\'ecessaire car $\tau$ n'est pas un morphisme de ${\cal D}_{X}[s_1, \ldots ,s_p]$-Module.\\
  
  \noindent  \underline{Rappels sur les anneaux r\'eguliers} : Soit $A$ un anneau r\'egulier (non n\'ecessairement commutatif)
  qui v\'erifie la condition d'Auslander : pour tout module \`a gauche de type fini, tout sous-module \`a droite de
  $Ext^k_A(E,A)$ est de grade sup\'erieur  ou \'egal \`a $k$. Notons que cette condition est automatiquement v\'erifi\'ee si $A$ est un anneau commutatif
  r\'egulier et qu'elle  est \`a la base des propri\'et\'es du nombre grade et de la filtration d'un module par le complexe bidualisant.
  Rappelons ces propri\'et\'es (voir \cite{Bj2} appendice IV)~:
  
 \begin{enumerate}
 \item Si $ 0 \rightarrow E' \rightarrow E \rightarrow E'' \rightarrow 0$ esf une suite exacte de $A$-module \`a gauche de type fini, 
 ${\rm grade}\, E  = {\rm inf} \, ({\rm grade}\, E' ,{\rm grade}\, E'' ) $.
 \item Soit $0\subset G_l\subset \cdots \subset G_0=E$ le filtration d'un $A$-module \`a gauche de type fini par le complexe bidualisant, il existe unes suite exacte fonctorielle en $E$~:
 $$  0 \longrightarrow \frac{G_k}{G_{k+1}}\longrightarrow Ext^k_A(Ext^k_A(E,A),A) \longrightarrow Q_k \longrightarrow 0 \; , $$ 
 o\`u ${\rm grade}\, Q_k\geq k+2$.
  \item Si ${\rm grade}\, E  = k$, alors ${\rm grade}\, Ext^k_A(E,A ) = k$.
  \item Si $Ext^v_A(Ext^v_A(E,A),A)\neq 0$, ce module est $A$-module pur de grade $v$.
  \item $E$ est pur de grade $k$ si et seulement si $0=G_l =\ldots = G_{k+1}$  et $G_k = \ldots = G_0=E$.
   Dans ce cas nous avons la suite exacte~:
 $$  0 \longrightarrow E \stackrel{i(E)}{\longrightarrow} Ext^k_A(Ext^k_A(E,A),A) \longrightarrow Q_k \longrightarrow 0 \; , $$ 
  o\`u ${\rm grade}\, Q_k\geq k+2$.
\end{enumerate}

\begin{definition} Soit $E$ un $A$-module de type fini pur de grade $k$, nous disons que $(\phi,E')$ est une extension pure docile de $E$
si  $\phi : E \rightarrow E'$ est un morphisme injectif de $A$-modules,
si $E'$ est  pur de ${\rm grade}\, E'  = k$ et ${\rm grade}\, (E'/\phi (E))  \geq k+2$.
  \end{definition}
  
La terminologie extension pure docile est une  traduction de tame pure extension ( voir \cite{Bj2} appendice IV).
Suivant les rappels, si $E$ est  pur de grade $k$, le morphisme naturel  $E \stackrel{i(E)}{\hookrightarrow} Ext^k_A(Ext^k_A(E,A),A)$ est une extension pure docile. En fait, cette extension naturelle est universelle au sens suivant~: 

 \begin{proposition} \label{pet}
 Soit $E$ un $A$-module de type fini pur de grade $k$ et  $(\phi,E')$   une extension pure docile de $E$. Alors, il existe un unique morphisme de $A$-module $\tilde {\phi} : E \hookrightarrow Ext^k_A(Ext^k_A(E,A),A)$ tel que le diagramme suivant commute~:
  $$
  \begin{array}{lcl} 
  E & \stackrel{\phi}{\hookrightarrow}   & E' \\ \downarrow i(E) & & \downarrow \tilde{\phi} \\  Ext^k_A(Ext^k_A(E,A),A)& = & Ext^k_A(Ext^k_A(E,A),A)
 \end{array} 
 $$ 
 De plus, $\tilde {\phi}$ est une extension pure docile.
  \end{proposition} 
  \noindent {\bf Preuve : } (proposition 2.9, \cite{Bj2} appendice IV)  Ce r\'esultat est la clef de la proposition \ref{psfp}. Donnons une preuve. Consid\'erons la suite exacte :
 $$   0 \rightarrow E \stackrel{\phi}{\rightarrow}  E'    \rightarrow  \frac{E'}{\phi (E)} \rightarrow 0\; . $$ 
 Nous obtenons vu les hypoth\`eses sur les nombres grade :
 $$ 0 = Ext^k_A(\frac{E'}{\phi (E)},A) \rightarrow  Ext^k_A(E',A)  \stackrel{Ext^k_A(\phi,A)}{\rightarrow}   Ext^k_A(E,A) \rightarrow
Ext^{k+1}_A(\frac{E'}{\phi (E)},A)=0$$ 
Donc, $Ext^k_A(\phi,A)$ est un isomorphisme. Il en r\'esulte le diagramme commutatif~:
$$
  \begin{array}{lcl} 
  E & \stackrel{\phi }{\longrightarrow}  & E' \\ \downarrow i(E) & & \downarrow i(E') \\  Ext^k_A(Ext^k_A(E,A),A)& \stackrel{ Ext^k_A(Ext^k_A(\phi,A),A)} \simeq & Ext^k_A(Ext^k_A(E',A),A)
 \end{array} 
 $$
 L'existence de $\tilde{\phi}$ s'en d\'eduit. L'application  $\tilde{\phi}$ est clairement une extension pure docile. Il reste \`a montrer l'unicit\'e de $\tilde{\phi}$. Si $\tilde{\phi}'$ est un deuxi\`eme morphisme 
 v\'erifiant $  \tilde{\phi}' \circ \phi = i(E)$, nous avons 
 $ ( \tilde{\phi}' -   \tilde{\phi})\circ \phi = 0$ et $\phi (E) \subset \ker ( \tilde{\phi}' -   \tilde{\phi})$. Comme $E'/ \ker ( \tilde{\phi}' -   \tilde{\phi})$
 s'injecte par $\tilde{\phi}' -   \tilde{\phi}$  dans $Ext^k_A(Ext^k_A(E',A),A)$, nous d\'eduisons de la puret\'e du double ext    que si   $E'/ \ker ( \tilde{\phi}' -   \tilde{\phi})$ est non nul, il est  de grade $k$. Mais $E'/\phi (E)$ est de grade sup\'erieur  ou \'egal \`a $k+2$ et se surjecte dans $E'/ \ker ( \tilde{\phi}' -   \tilde{\phi})$. Il en r\'esulte $\ker ( \tilde{\phi}' -   \tilde{\phi}) = E$ et $\tilde{\phi}' = \tilde{\phi}$.\\

  \noindent {\bf Preuve de la proposition \ref{psfp}  : }  Nous d\'etaillons  la preuve   au niveau des fibres de  $M$ en un point $x_0$ de $X$. La coh\'erence de $\tilde{L} $ r\'esultera du fait que  les sous-modules de ma famille
   ${\cal G}$ peuvent \^etre consid\'er\'es en utilisant une version faisceautique de la proposition \ref{pet} comme des sous-modules de  
  $$  Ext^{ {\rm dim}\, X + 1 }_{{\cal D}_{X }[s_1, \ldots ,s_p]}(Ext^{{\rm dim}\, X + 1 }_{{\cal D}_{X }[s_1, \ldots ,s_p]}(   {\cal D}_{X }[s_1, \ldots ,s_p] m f_1^{s_1 } \ldots  f_p^{s_p }  ,{\cal D}_{X}[s_1, \ldots ,s_p]) \quad .$$

  \noindent  {\bf Si }
  $$  \frac {     {\cal D}_{X,x_0}[s_1, \ldots ,s_p] m f_1^{s_1-1} \ldots  f_p^{s_p-1}   }{   {\cal D}_{X,x_0}[s_1, \ldots ,s_p] m f_1^{s_1 } \ldots  f_p^{s_p }  }  
   $$
   {\bf est pur  de grade $k+1$},  nous v\'erifions que   $\tilde{L}_{x_0}=  {\cal D}_{X,x_0}[s_1, \ldots ,s_p] m f_1^{s_1   } \ldots  f_p^{s_p  }  $. Nous notons en utilisant l'action de $\tau$ que  les quotients :
    $$  \frac {     {\cal D}_{X,x_0}[s_1, \ldots ,s_p] m f_1^{s_1+k} \ldots  f_p^{s_p+k}   }{   {\cal D}_{X,x_0}[s_1, \ldots ,s_p] m f_1^{s_1+k+1} \ldots  f_p^{s_p+k+1}  }  
   $$ 
   ont m\^eme nombre grade. \\
   
    \noindent   {\bf Sinon}, nous avons vu que : $$  \frac {     {\cal D}_{X,x_0}[s_1, \ldots ,s_p] m f_1^{s_1-1} \ldots  f_p^{s_p-1}   }{   {\cal D}_{X,x_0}[s_1, \ldots ,s_p] m f_1^{s_1 } \ldots  f_p^{s_p }  }  
   $$
   est de grade sup\'erieur  ou \'egal \`a $k+1$. Ce quotient admettrait donc un sous-module de grade sup\'erieur  ou \'egal \`a $k+2$.
   Il existerait ainsi $L$ un $ {\cal D}_{X,x_0}[s_1, \ldots ,s_p] $-module tel que~: 
  $$ {\cal D}_{X,x_0}[s_1, \ldots ,s_p] m f_1^{s_1 } \ldots  f_p^{s_p } \stackrel{\neq}{\subset}  L \subset
 {\cal D}_{X,x_0}[s_1, \ldots ,s_p] m f_1^{s_1 -1} \ldots  f_p^{s_p - 1} $$
 et
$$  {\rm grade} \, \frac{L}{ {\cal D}_{X,x_0}[s_1, \ldots ,s_p] m f_1^{s_1 } \ldots  f_p^{s_p } } \geq k+2 \; .$$ 
   Consid\'erons alors la famille ${\cal G}_{x_0}$  des   ${\cal D}_{X,x_0}[s_1, \ldots ,s_p]$ modules $L $  de type fini v\'erifiant
$$ {\cal D}_{X,x_0}[s_1, \ldots ,s_p] m f_1^{s_1 } \ldots  f_p^{s_p } \stackrel{\neq}{\subset}  L \subset
M_{x_0}[\frac{1}{F},s_1, \ldots ,s_p] f_1^{s_1 } \ldots  f_p^{s_p } 
 $$
 et
$$  {\rm grade} \, \frac{L}{ {\cal D}_{X,x_0}[s_1, \ldots ,s_p] m f_1^{s_1 } \ldots  f_p^{s_p } } \geq k+2 \; .$$ 
   Le module $M_{x_0}[\frac{1}{F},s_1, \ldots ,s_p] f_1^{s_1 } \ldots  f_p^{s_p } $ 
   n'est pas de type fini. Mais pour tout sous-module $N$ de type fini, il existe un entier $r$  tel que $N$  est un sous-module
de   $ {\cal D}_{X,x_0}[s_1, \ldots ,s_p] m f_1^{s_1 -r} \ldots  f_p^{s_p - r}   $ qui est un module  pur de grade ${\rm dim}\, X$ 
(voir proposition \ref{ppdmfs}). Ainsi, tous les sous-modules de type fini de $M_{x_0}[\frac{1}{F},s_1, \ldots ,s_p] f_1^{s_1 } \ldots  f_p^{s_p } $
sont purs de   grade ${\rm dim}\, X$ .
\begin{lemme}  ${\cal G}_{x_0}$ admet un unique plus grand \'el\'ement.
\end{lemme}

  \noindent {\bf Preuve du lemme  : }  (\cite{Bj2}, proposition 2.10 appendice IV) Rappelons cette preuve qui repose sur le probl\`eme 
   universel des extensions pures dociles. Par les propri\'et\'es des nombres grade, la somme de deux \'el\'ements de 
    ${\cal G}_{x_0}$ est dans  ${\cal G}_{x_0}$. Il suffit donc de montrer que toute suite croissante  d'\'el\'ements de  ${\cal G}_{x_0}$ stationne. Soit $L_v$
    une telle suite, notons $j_v$ et $j_{v,v+1}$ les inclusions :
    $$ E = {\cal D}_{X,x_0}[s_1, \ldots ,s_p] m f_1^{s_1 } \ldots  f_p^{s_p }  \stackrel{j_v }{\longrightarrow} L_v \stackrel{j_{v,v+1} }{\longrightarrow}
    L_{v+1} \; .$$ 
    Suivant le probl\`eme universel des extensions pures dociles, nous avons le diagramme commutatif~:
  $$
  \begin{array}{lcr} 
  E = {\cal D}_{X,x_0}[s_1, \ldots ,s_p] m f_1^{s_1 } \ldots  f_p^{s_p } & \stackrel{j_v }{\longrightarrow}  & L_v   \\
  \downarrow i(E) & & \swarrow \tilde{j_v}  \\
  Ext^{ {\rm dim}\, X + 1 }(Ext^{{\rm dim}\, X + 1 }(E ,{\cal D}_{X,x_0}[s_1, \ldots ,s_p]),{\cal D}_{X,x_0}[s_1, \ldots ,s_p])&  &
       \end{array} 
 $$
  
 Par unicit\'e de $\tilde{j_v}$, nous obtenons $ \tilde{j_v} = \tilde{j}_{v+1} \circ j_{v,v+1} $. Ainsi, les images des morphismes $ \tilde{j_v}$
 forment une suite croissante de sous-modules d'un module noeth\'erien de type fini.  Cette suite d'images et donc les $L_v$ stationnent.\\
 
\noindent {\bf Notation :} Nous notons $\tilde {L}_{x_0}$ le plus grand \'el\'ement de ${\cal G}_{x_0}$.\\
  
  Par finitude de $\tilde {L}_{x_0}$ , il existe un entier $r $ tel que :
  $$\tilde {L}_{x_0} \subset  {\cal D}_{X,x_0}[s_1, \ldots ,s_p] m f_1^{s_1 -r} \ldots  f_p^{s_p - r}  $$  
Nous allons montrer que $\tilde {L}_{x_0} $ satisfait les propri\'et\'es de la proposition  \ref{psfp}.\\

Posons pour simplifier pout tout entier $r$ et $s= (s_1,\ldots ,s_p)$~:
$$ E(s+r) = {\cal D}_{X,x_0}[s_1, \ldots ,s_p] m f_1^{s_1 +r} \ldots  f_p^{s_p +r}\; .$$

  \noindent { \bf Preuve  de $\tau ( \tilde {L}_{x_0}) \subset \tilde {L}_{x_0}$ : }   Consid\'erons le morphisme surjectif~:
  $$\frac{\tau ( \tilde {L}_{x_0} ) }{E(s+1)} \longrightarrow \frac{\tau ( \tilde {L}_{x_0} ) + E(s) }{E(s)} \; . $$
Mais :
$$ \frac{\tau ( \tilde {L}_{x_0} ) }{E(s+1)} =  \tau (   \frac{  \tilde {L}_{x_0}  }{E( s)} ) $$
(nous passons au quotient l'action de $\tau$) et donc~:
$$ {\rm grade}\, \frac{\tau ( \tilde {L}_{x_0} ) + E(s) }{E(s)} \geq  {\rm grade}\, \frac{\tau ( \tilde {L}_{x_0} ) }{E(s+1)}  =  {\rm grade}\, \frac{  \tilde {L} _{x_0} }{E( s)} \geq  {\rm dim}\, X +2 \; .$$
Ainsi, $\tau ( \tilde {L}_{x_0} ) + E(s) \in {\cal G}  $,   $\tau ( \tilde {L}_{x_0} ) + E(s) \subset \tilde {L} _{x_0} $ et  $\tau ( \tilde {L}_{x_0} )  \subset \tilde {L} _{x_0} $.\\

  \noindent { \bf Preuve  de $ \displaystyle {\rm grade}\, \frac{  \tilde {L}_{x_0}   }{ \tau ( \tilde {L} _{x_0} )} =  {\rm grade}\, \frac{  E(s) }{E(s+1) }  \geq    {\rm dim}\, X +1$ : }  
Consid\'erons les  inclusions $ E(s)\subset \tilde {L}_{x_0} \subset E(s-r)$ et $ \tau ( \tilde {L}_{x_0} ) \subset \tilde {L}_{x_0} $~:
$$ \frac{  \tau ( \tilde {L}_{x_0}) }{E(s+1) }  \hookrightarrow \frac{   \tilde {L} _{x_0} }{E(s+1) }  \hookrightarrow \frac{  E(s-r)  }{E(s+1) } \; , $$
$$ \frac{    E(s+1)  }{ \tau ^{r+1 }( \tilde {L}_{x_0} ) }  \hookrightarrow \frac{  E(s )   }{ \tau ^{r+1 }( \tilde {L}_{x_0} ) }  
\hookrightarrow \frac{   \tilde {L}_{x_0}     }{ \tau ^{r+1 }( \tilde {L} _{x_0}) }\; . $$ 
Donc, $ \tilde {L}_{x_0}  /   \tau ( \tilde {L} _{x_0})  $ (resp. $E(s)/E(s+1) $) est un sous-facteur de $E(s-r)/ E(s+1)$ (resp. $ \tilde {L} _{x_0}/ \tau ^{r+1 }( \tilde {L}_{x_0} )$. Nous en d\'eduisons~:
$${\rm grade}\,   \frac{\tilde {L}_{x_0}  }{  \tau ( \tilde {L} _{x_0})} \geq {\rm grade}\,  \frac{E(s-r)}{E(s+1)}   \; {\rm et} \;  
{\rm grade}\,  \frac{E(s)}{E(s+1) } \geq {\rm grade}\,  \frac{\tilde {L}_{x_0}} { \tau ^{r+1 }( \tilde {L} _{x_0}) }\; .$$
Par action de $\tau$,  pour  $k$ entier, les grades de $E(s+k)/E(s+k+1)$ et $ \tau ^{k}  ( \tilde {L}_{x_0}  ) / \tau ^{k+1}  ( \tilde {L}_{x_0} )$ sont  ind\'ependants de $k$. Nous en d\'eduisons par r\'ecurrence que $E(s-r)/ E(s+1)$ (resp.   $ \tilde {L}_{x_0} / \tau ^{r+1 }( \tilde {L}_{x_0})$  est de m\^eme grade que
$  {E(s)}/{E(s+1) }$  (resp. $\tilde {L}_{x_0}  /   \tau ( \tilde {L}_{x_0} )$. Il en r\'esulte~:
$$  {\rm grade}\,   \frac{\tilde {L}_{x_0}  }{  \tau ( \tilde {L}_{x_0} )} = {\rm grade}\,  \frac{E(s)}{E(s+1) } \; .$$ 
Enfin, nous avons montr\'e que $ \displaystyle {\rm grade}\,  \frac{E(s)}{E(s+1) } \geq    {\rm dim}\, X +1$ (voir lemme \ref{lmgef}).\\

\noindent { \bf Preuve  de $   \displaystyle  \frac{  \tilde {L}_{x_0}   }{ \tau ( \tilde {L}_{x_0}  )} $ pur de grade $   {\rm dim}\, X +1$ : }  
Consid\'erons un ${\cal D}_{X,x_0}[s_1, \ldots ,s_p] $-module $L'$ tel que
$$ E(s+1 ) \subset \tau ( \tilde {L}_{x_0}  ) \stackrel{\neq}{\subset} L' \subset \tilde {L}_{x_0} \; . $$ 
Nous avons la suite exacte~:
$$ 0 \longrightarrow \frac{ \tau ( \tilde {L}_{x_0}  )  }{E(s+1 ) }
\longrightarrow  \frac{   L'    }{E(s+1 ) }
\longrightarrow  \frac{   L'    }{\tau ( \tilde {L}_{x_0} )}
\longrightarrow  0 \; .$$ 
Par maximalit\'e de $\tilde {L}_{x_0} $~:
$$ {\rm grade}\,     \frac{   L'    }{E(s+1 ) } =  {\rm grade}\,     \frac{  \tau ^{-1} (L' )   }{E(s ) } \leq  {\rm dim}\, X +1\; .$$
Or,
$$  {\rm grade}\,   \frac{ \tau ( \tilde {L}_{x_0}  )  }{E(s+1 ) } = {\rm grade}\,  \frac{   \tilde {L}_{x_0}     }{E(s   ) } \geq  {\rm dim}\, X + 2\; .$$
Il en r\'esulte~:
$$   {\rm grade}\,  \frac{   L'    }{\tau ( \tilde {L}_{x_0} )} \leq {\rm dim}\, X +1$$
et en prenant $  L'   =  \tilde {L}_{x_0} $~:
$$ {\rm grade}\,  \frac{ \tilde {L}_{x_0}   }{\tau ( \tilde {L}_{x_0} )} \leq {\rm dim}\, X +1 $$
Mais 
$$   {\rm grade}\,  \frac{   L'    }{\tau ( \tilde {L}_{x_0} )} \geq {\rm grade}\,  \frac{ \tilde {L}_{x_0}   }{\tau ( \tilde {L} _{x_0})} \geq {\rm dim}\, X +1 \; .
$$  
Ainsi,
$$  {\rm grade}\,  \frac{   L'    }{\tau ( \tilde {L}_{x_0} )} ={\rm grade}\,  \frac{ \tilde {L}_{x_0}   }{\tau ( \tilde {L}_{x_0} )} = {\rm dim}\, X +1 \; .  $$
Tout sous-module non nul de $\tilde {L}_{x_0}   / \tau ( \tilde {L}_{x_0} )$ est donc de grade ${\rm dim}\, X+1$ et  $\tilde {L}_{x_0}   / \tau ( \tilde {L}_{x_0} )$ est donc pur de grade 
${\rm dim}\, X+1$.\\

 \begin{remarque} \label{rp} Notons que $ \tilde {L} = {\cal D}_{X }[s_1, \ldots ,s_p] m f_1^{s_1 } \ldots  f_p^{s_p } $ si et seulement si les 
les fibres  de $$  \frac {     {\cal D}_{X }[s_1, \ldots ,s_p] m f_1^{s_1 } \ldots  f_p^{s_p }   }{   {\cal D}_{X }[s_1, \ldots ,s_p] m f_1^{s_1 + 1 } \ldots  f_p^{s_p + 1}  }  
   $$
   sont des modules  purs de grade ${\rm dim}\, X +1$. 
   Cette  condition \'equivaut encore  :
   $$  Ext^{i}_{{\cal D}_{X }[s_1, \ldots ,s_p]}(Ext^{i }_{{\cal D}_{X }[s_1, \ldots ,s_p]}(   \frac {     {\cal D}_{X }[s_1, \ldots ,s_p] m f_1^{s_1 } \ldots  f_p^{s_p }   }{   {\cal D}_{X }[s_1, \ldots ,s_p] m f_1^{s_1 + 1 }
   \ldots  f_p^{s_p + 1}  }    ,{\cal D}_{X}[s_1, \ldots ,s_p])) \neq 0  \Longrightarrow  \; i = {\rm dim}\, X +1 \; .$$
    \end{remarque}

 \noindent { \bf Preuve :}  Voir la d\'efinition de  $\tilde {L} $.\\

Donnons maintenant quelques propri\'et\'es de $\tilde {L} $.

\begin{proposition}\label{pplt} Le  ${\cal D}_{X }[s_1, \ldots ,s_p] $-module  $\tilde {L} $ construit dans la proposition pr\'ec\'edente v\'erifie~:
\begin{enumerate}
\item $\displaystyle {\rm car}^{\sharp}\,  \frac{ \tilde {L}   }{\tau ( \tilde {L} )} = {\rm car}^{\sharp}\, \frac {     {\cal D}_X[s_1, \ldots ,s_p] m f_1^{s_1} \ldots  f_p^{s_p}   }{   {\cal D}_X[s_1, \ldots ,s_p] m f_1^{s_1+1} \ldots  f_p^{s_p+1}  }$,
\item localement au voisinage de tout point $x_0$ de $X$, il existe un entier $r$ tel que :
$$\displaystyle {\rm car}^{\rm rel}\,  \frac{ \tilde {L}   }{\tau ( \tilde {L} )} \subset \bigcup_{k=0}^r {\rm car}^{\rm rel }\, \tau ^k \left( \frac {     {\cal D}_X[s_1, \ldots ,s_p] m f_1^{s_1} \ldots  f_p^{s_p}   }{   {\cal D}_X[s_1, \ldots ,s_p] m f_1^{s_1+1} \ldots  f_p^{s_p+1}  } \right)\; ,$$
\item   localement au voisinage de tout point $x_0$ de $X$, il existe un entier $r$ tel que :
$$\displaystyle{\rm car}^{\rm rel}\,  \frac {     {\cal D}_X[s_1, \ldots ,s_p] m f_1^{s_1} \ldots  f_p^{s_p}   }{   {\cal D}_X[s_1, \ldots ,s_p] m f_1^{s_1+1} \ldots  f_p^{s_p+1}  } \subset 
 \bigcup_{k=0}^r {\rm car}^{\rm rel }\, \tau ^{-k }  \left(     {\rm car}^{\rm rel}\,  \frac{ \tilde {L}   }{\tau ( \tilde {L} )}   \right)\; .$$

\end{enumerate}
\end{proposition}

  \noindent { \bf Preuve  :} Localement au voisinage de tout point $x_0$ de $X$, il existe un entier $r$ tel que 
  $\tilde {L}  / \tau ( \tilde {L} )$ (resp. $E(s)/E(s+1)$ est quotient de sous-modules de $E(s-r)/E(s+1)$ (resp. $   \tilde {L}  / \tau  ^{r+1}( \tilde {L} )$). Nous utilisons de plus que pour $L$ sous-module coh\'erent de
  $M [1/F,s_1, \ldots ,s_p]   f_1^{s_1} \ldots  f_p^{s_p}$~:
  $$ {\rm car}^{\sharp}\,  \frac{  {L}   }{\tau (  {L} )} = {\rm car}^{\sharp}\, \frac{  {\tau ^k(   {L} )}   }{\tau ^{k+1}(   {L} )}
\quad {\rm et} \quad {\rm car}^{\rm rel}\, \frac{  {\tau ^k(   {L} )}   }{\tau ^{k+1}(   {L} )} =  \tau ^{-k} \left(     {\rm car}^{\rm rel}\,  \frac{  L   }{\tau (  L  )}   \right)\; .$$

Rappelons quelques propri\'et\'es des modules purs sur un anneau $A$ filtr\'e positivement dont le gradu\'e ${\rm gr} \, A$ est commutatif r\'egulier. Si un  ${\rm gr} \, A$-module
de type fini est pur, ses id\'eaux associ\'es co\"{\i}ncident avec les id\'eaux premiers minimaux de son support et la hauteur de ces id\'eaux est constante \'egale au grade du module (voir 
\cite{Bj2}, appendice IV proposition 3.7). Cette condition caract\'erise en fait les modules purs (voir  \cite{Bj2}, appendice IV remarque 3.8). Si $N$  est un $A$-module pur, il existe une bonne filtration de $N$ tel que le gradu\'e de $N$ soit un ${\rm gr} \, A$-module pur (voir  \cite{Bj2}, appendice IV theorem 4.11). La cons\'equence est que ${\cal J}\,(N)$ la racine de l'annulateur de
${\rm gr} \, N$ a tous ses id\'eaux associ\'es de m\^eme hauteur. Prendre garde que  pour $p$ id\'eal premier de ${\rm gr} \, A$, nous n'avons pas n\'ecessairement 
${\rm dim }\, ( {\rm gr} \, A / p ) + {\rm ht }\, p = {\rm dim} \,  {\rm gr} \, A $.\\

En particulier si $N$ est un ${\cal D}_{X,x_0}[s_1, \ldots ,s_p] $-module pur de grade ${\rm dim}\, X +1$, il existe une bonne filtration di\`ese de $N$
tel que $ {\rm gr}^{\sharp } \, N$ soit pur de m\^eme grade. Les id\'eaux pemiers associ\'es de $ {\rm gr}^{\sharp } \, N  $ sont homog\`enes en $(\xi,s)$ et leurs hauteurs  co\"{\i}ncident avec la codimension des germes des espaces analytiques qu'ils d\'efinissent sur $T^{\ast}X \times {\bf C}^p $.  
De la proposition \ref{pplt}, nous d\'eduisons en prenant $N =  \tilde {L}  / \tau ( \tilde {L} )$~:

\begin{corollaire} Les composantes irr\'eductibles de   
$${\rm car}^{\sharp}\,  \frac{ \tilde {L}   }{\tau ( \tilde {L} )} =   {\rm car}^{\sharp}\,   \frac{{\cal D}_X[s_1, \ldots ,s_p] m f_1^{s_1} \ldots  f_p^{s_p}}{{\cal D}_X[s_1, \ldots ,s_p] m f_1^{s_1+1} \ldots  f_p^{s_p+1}} \; .$$
sont toutes de dimension ${\rm dim}\, X +p -1  $.  
\end{corollaire}

\begin{corollaire}\label{cvcl} Il existe une famille $(X_{\alpha})_{\alpha \in A }$    de sous-espaces analytiques  de $X$
et une famille  $(S'_{\alpha})_{\alpha \in A }$ de  r\'eunion d'hypersurfaces alg\'ebriques de   ${\bf C}^p$
telles que :
$$  {\rm car}^{\rm rel }\,  (  \frac{ \tilde {L}   }{\tau ( \tilde {L} )} )  =  \cup_{\alpha \in A'}T^{\ast}_{Y_{\alpha}}X \times  S'_{\alpha} \; .$$
\end{corollaire}

\noindent { \bf Preuve :}   Le  ${\cal D}_X[s_1, \ldots ,s_p] $-Module  $ \tilde {L}  / \tau ( \tilde {L} )  $ \'etant major\'e par une lagrangienne, il existe suivant la proposition   \ref{psvcrmml}  une famille $(Y_{\alpha})_{\alpha \in A'}$ (resp. $S'_{\alpha}$)   de sous-espaces analytiques (resp. vari\'et\'es alg\'ebriques de ${\bf C}^p$)
telle que $  {\rm car}^{\rm rel }\,  ( \tilde {L}  / \tau ( \tilde {L} ) )  =  \cup_{\alpha \in A'}T^{\ast}_{Y_{\alpha}}X \times  S'_{\alpha}$ .
Alors, pour tout $\alpha$,
les $S'_{\alpha}$ sont de dimension $p-1$. Nous consid\'erons le module pur $N =  \tilde {L}  / \tau ( \tilde {L} )$.  Il admet donc une bonne filtration relative tel que
   $ {\rm gr}^{\rm rel } \, N$ soit pur. Si $p$ est un id\'eal premier minimal du support de ce gradu\'e, le quotient   ${\rm gr}^{ \rm rel } \,( {\cal D}_{X,x_0}[s_1, \ldots ,s_p] )  / p  $ 
est pur  de grade  ${\rm dim}\, X +1$. Il en r\'esulte que le  sous-espaces analytique  qu'il  d\'efinit est de dimension  ${\rm dim}\, X +  p-1$. Toutes les composantes de 
  $ {\rm car}^{\rm rel}\, N$ sont donc de dimension ${\rm dim}\, X + p-1$. Les   $S'_{\alpha}$ sont donc  de dimension $p-1$.

 \begin{remarque} Nous verrons dans le paragraphe suivant que les  $S'_{\alpha}$ sont en fait des hyperplans affines.
  Il en r\'esultera que   la vari\'et\'e des z\'eros de l'id\'eal de Bernstein
${\cal B}_{x_0}( \tilde {L}  / \tau ( \tilde {L} )) $ est une  r\'eunion d'hyperplans affines. 

 \end{remarque} 
  
  \subsection{Remarques sur l'id\'eal de Bersntein ${\cal B}(m,x_0, f_1, \ldots ,f_p)$}
  
  \subsubsection{Cas o\`u $m$ est une section d'un Module holonome}
  
  $M$ d\'esigne toujours  un ${\cal D}_X$-Module holonome engendr\'e par une section $m$ et   $f_1, \ldots ,f_p$ sont  des fonctions analytiques sur $X$.\\
  
  Rappelons, voir d\'efinition \ref{dfnibm}, que ${\cal B}(m,x_0, f_1, \ldots ,f_p)$ d\'esigne l'id\'eal    des polyn\^omes   $b$ de   ${\bf C}[s_1, \ldots ,s_p]$ 
  v\'erifiant  au voisinage de $x_0$ : 
$$ (\ast) \; \quad \;   b (s_1, \ldots ,s_p) m f_1^{s_1} \ldots  f_p^{s_p} \in {\cal D}_X[s_1, \ldots ,s_p] m f_1^{s_1+1} \ldots  f_p^{s_p+1} \; .$$
Ces polyn\^omes sont appel\'es  polyn\^omes de Bernstein de $(m, f_1, \ldots ,f_p) $ au voisinage de $x_0$.\\

    
  \begin{proposition} \label{pcspfl} (C. Sabbah proposition 1.2 de \cite{S2}) Pour toute section $m$ d'un Module holonome,  pour tout $x_0 \in X$, il existe  un nombre fini de formes lin\'eaires ${\cal H}$ \`a coefficients premiers entre eux dans ${\bf N}$ telles que :
    $$\prod_{H\in {\cal H} } \prod_{i\in I_{\cal H} } (H(s) + \alpha _{H,i}) \in {\cal B}(m,x_0, f_1, \ldots ,f_p) $$
    o\`u $ \alpha _{H,i}$ sont des nombres complexes. 
     \end{proposition}
     
      La preuve de C. Sabbah  repose sur la proposition 2.2.3  de \cite{S1},   voir \'egalement \cite{S3}. 
    Si $M = {\cal O}_X$, comme indiqu\'e dans \cite{S3},  les complexes peuvent \^etre pris rationels positifs (  voir   une preuve dans l'article \cite{Go} de A. Gyoja).\\

\begin{proposition} \label{psib}
Il existe 
      $(X_{\alpha})_{\alpha \in A}$ (resp. $S_{\alpha}$) une  famille finie de sous-espaces analytiques (resp. vari\'et\'es alg\'ebriques de ${\bf C}^p$
      telles que 
$$  {\rm car}^{\rm rel }\,  \left(  \frac{{\cal D}_X[s_1, \ldots ,s_p] m f_1^{s_1} \ldots  f_p^{s_p}}{{\cal D}_X[s_1, \ldots ,s_p] m f_1^{s_1+1} \ldots  f_p^{s_p+1}}  \right)  =  \cup_{\alpha \in A}T^{\ast}_{X_{\alpha}}X \times  S_{\alpha}   \; .$$
 avec les conditions :
 
 \begin{itemize}
 \item chaque vari\'et\'e alg\'ebrique  $S_{\alpha}$ est de dimension $p-1$,  
 \item les composantes irr\'eductibles de 
 dimension $p-1$ de chaque $S_{\alpha}$  sont des hyperplans affines  $H_{\alpha,\beta}$  de directions les noyaux de  formes lin\'eaires \`a coefficients premiers entre eux dans ${\bf N}$,
 \item  toute   composante  irr\'eductible   des  $S_{\alpha}$  de    dimension strictement inf\'erieure \`a $p-1$ est  contenue  dans un hyperplan affine $\tau ^k ( H_{\alpha,\beta}  )$ o\`u $k \in {\bf Z}$ et 
  $\tau$ la translation $(s_1, \ldots, s_p) \mapsto (s_1+1, \ldots, s_p+1)$. \\
  \end{itemize} 
 \end{proposition}
   
   \noindent { \bf Preuve :}  
    Consid\'erons le ${\cal D}_X[s_1, \ldots ,s_p]$-module $\tilde{L}$ construit \`a la proposition \ref{psfp}. Au voisinage de tout point $x_0$, il existe un entier $r$ tel que ~:
     $$ {\cal D}_X[s_1, \ldots ,s_p] m f_1^{s_1} \ldots  f_p^{s_p} \subset \tilde{L}  \subset 
     {\cal D}_X[s_1, \ldots ,s_p] m f_1^{s_1-r} \ldots  f_p^{s_p-r }  \; .$$  
     
Suivant le corollaire \ref{cvcl}, il existe $(X'_{\alpha})_{\alpha \in A }$ (resp. $S'_{\alpha}$) une famille de sous-espaces analytiques (resp. de vari\'et\'es alg\'ebriques de ${\bf C}^p$ de dimension $p-1$)  tels que 
$$  {\rm car}^{\rm rel }\,  ( \tilde {L}  / \tau ( \tilde {L} ) )  =  \cup_{\alpha \in A'}T^{\ast}_{Y_{\alpha}}X \times  S'_{\alpha}\; .$$
Consid\'erons $b_S$ un polyn\^ome de $ {\cal B}(m,x_0, f_1, \ldots ,f_p) $ parmi ceux dont l'existence est assur\'ee dans la proposition \ref{pcspfl}.
Comme    $b_S(s  -r)\cdots b_S(s  ) \tilde{L}  \subset \tau (   \tilde{L})$, ce produit  $b_S(s  -r)\cdots b_S(s  ) $  s'annule   sur les $S'_{\alpha}$
 tels que   $x_0 \in X_{ \alpha}$. Vu la dimension des 
  $S'_{\alpha}$, cela montre que les  $S'_{\alpha}$ sont des r\'eunions d'hyperplans affines. \\
  
  D'autre part  suivant la proposition \ref{pefml},   il existe une vari\'et\'e lagrangienne $\Lambda$ de $T^{\ast}X$ telle que ~:
$$ {\rm car}^{\rm rel} \,  {\cal D}_X[s] m f_1^{s_1} \ldots  f_p^{s_p} = \Lambda \times {\bf C}^p \; .$$
Il existe donc     $(X_{\alpha})_{\alpha \in A}$ (resp. $S_{\alpha}$) une famille de sous-espaces analytiques (resp. de vari\'et\'es alg\'ebriques de ${\bf C}^p$) telles que 
$$  {\rm car}^{\rm rel }\,  \left(  \frac{{\cal D}_X[s_1, \ldots ,s_p] m f_1^{s_1} \ldots  f_p^{s_p}}{{\cal D}_X[s_1, \ldots ,s_p] m f_1^{s_1+1} \ldots  f_p^{s_p+1}}  \right)  =  \cup_{\alpha \in A}T^{\ast}_{X_{\alpha}}X \times  S_{\alpha}   \; .$$
Nous d\'eduisons alors de la proposition \ref{pplt} que $A'= A$ et $X_{\alpha}= Y_{\alpha}$ et ~:
$$ \bigcup  S'_{\alpha} \subset  \bigcup_{k=0}^r \tau ^k (  \bigcup  S_{\alpha} ) \;  {\rm et} \; 
  \bigcup   S_{\alpha}  \subset \bigcup_{k=0}^r \tau ^{-k} ( \bigcup  S'_{\alpha} ) \; .$$
  La proposition en r\'esulte.\\

\begin{definition}
Pour tout $x_0 \in X$, notons ${\cal H}(x_0,m)$   l'ensemble des directions des hyperplans $H_{\alpha,\beta}$ pour les  $\alpha$   tels que $x_0 \in X_{\alpha}$. Nous 
 appelons ces directions les  pentes de $(m,f_1, \dots,f_p)$ au voisinage de $x_0$.
 \end{definition}

  Suivant la proposition \ref{pvcrib},   $\bigcup_{x_0 \in X_{\alpha}} S_{\alpha} $
    est la vari\'et\'e des z\'eros de l'id\'eal de Bernstein de $  (m,  f_1, \ldots ,f_p)$. Nous en d\'eduisons :
     
  \begin{corollaire}\label{cpbpfl}
    Si un produit de  formes lin\'eaires affines appartient \`a $  {\cal B}(m,x_0, f_1, \ldots ,f_p)$, tout hyperplan vectoriel  de  ${\cal H}(x_0,m)$ est direction de l'un des facteurs.
De plus, il existe dans $  {\cal B}(m,x_0, f_1, \ldots ,f_p)$ un produit de formes lin\'eaires affines dont les directions sont exactement l'ensemble ${\cal H}(x_0,m)$ des 
    pentes de  $(m,f_1, \dots,f_p)$ au voisinage de $x_0$.
  \end{corollaire} 

 
 \noindent { \bf Preuve :}   Pour le premier point, il suffit d'observer que
 si $b \in {\cal B}(m,x_0, f_1, \ldots ,f_p)$, il s'annule sur les $S_{\alpha}$  intervenant dans la proposition \ref{psib}, donc sur ceux de dimension $p-1$ dont     les directions 
sont les directions des hyperplans
 de  $ {\cal H}$. 
  Pour le deuxi\`eme  point, il suffit de prendre un produit de formes lin\'eaires nulles sur les $S_{\alpha}$  
 intervenant dans la proposition \ref{psib} et d'utiliser la proposition \ref{pvcrib}.

  \begin{corollaire} 
  La vari\'et\'e  des z\'eros de ${\bf in }\,  {\cal B}(m,x_0, f_1, \ldots ,f_p)$  l'id\'eal   engendr\'e par les  parties homog\'enes de plus haut degr\'e  
     des \'el\'ements de $ {\cal B}(m,x_0, f_1, \ldots ,f_p)$ est la r\'eunion  des   pentes de  $(m,f_1, \dots,f_p)$ au voisinage de $x_0$..
      La  racine de ${\bf in }\,  {\cal B}(m,x_0, f_1, \ldots ,f_p)$  est en particulier un    id\'eal principal.
        \end{corollaire}

 
     \noindent { \bf Preuve :}    R\'esulte du corollaire  \ref{cpbpfl}.

     \begin{proposition} Soit $x_0 \in X$. Si le ${\cal D}_{X ,x_0}[s_1, \ldots ,s_p]$-module :
     $$  \frac {     {\cal D}_{X ,x_0}[s_1, \ldots ,s_p] m f_1^{s_1 } \ldots  f_p^{s_p }   }{   {\cal D}_{X,x_0 }[s_1, \ldots ,s_p] m f_1^{s_1 + 1 } \ldots  f_p^{s_p + 1}  }  
    $$
    est un   module   pur  de grade ${\rm dim}\, X +1$, la racine de l'id\'eal de Bernstien $  {\cal B}(m,x_0, f_1, \ldots ,f_p)$ est principal. C'est notamment le cas sous la condition :
    $$  Ext^{i}_{{\cal D}_{X,x_0 }[s_1, \ldots ,s_p]}(Ext^{i }_{{\cal D}_{X,x_0 }[s_1, \ldots ,s_p]}(   \frac {     {\cal D}_{X,x_0 }[s_1, \ldots ,s_p] m f_1^{s_1 } \ldots  f_p^{s_p }   }{   {\cal D}_{X,x_0 }[s_1, \ldots ,s_p] m f_1^{s_1 + 1 }
   \ldots  f_p^{s_p + 1}  }    ,{\cal D}_{X,x_0}[s_1, \ldots ,s_p])) \neq 0  \Longrightarrow  \; i = {\rm dim}\, X +1 \; .$$
  \end{proposition}
  
  \noindent { \bf Preuve :} Cela r\'esulte de la remarque \ref{rp} et du corollaire \ref{cvcl}. La deuxi\'eme partie de la proposition est un crit\`ere 
  de puret\'e (voir \cite{Bj2}, Appendice IV,  proposition 2.6). \\

     Pour terminer ce paragraphe, consid\'erons l'application :
     $$ {\rm exp}^{ 2i \pi . } : {\bf C}^p \longrightarrow  ({\bf C}^{\ast})^p \quad , \quad (s_1, \ldots ,s_p) \longmapsto (e^{2i\pi s_1}, \ldots ,e^{2i\pi s_p})\; .$$

   \begin{corollaire} L'image par l'application ${\rm exp}^{ 2i \pi . }$ de la  vari\'et\'es des z\'eros de $  {\cal B}(m,x_0, f_1, \ldots ,f_p)$ est une r\'eunion
        de sous-ensembles de $ ({\bf C}^{\ast})^p$  o\`u  chaque sous-ensemble est d\'efini  par une \'equation du  type :
  $$ (\sigma _1)^{a_1} \cdots (\sigma _p)^{a_p} = \alpha $$ 
  o\` u $(a_1, \ldots ,a_p)$ est une famille d'\'el\'ements de ${\bf N}$ premier entre eux et $\alpha $ un nombre complexe. L'ensemble des $(a_1, \ldots ,a_p)$ 
   est l'ensemble des coefficients d'\'equations  des   pentes de  $(m,f_1, \dots,f_p)$ au voisinage de $x_0$.  
 \end{corollaire}    
     
     \noindent { \bf Preuve :} Cela se d\'eduit de la proposition   \ref{psib}. En effet, si   $a_1s_1 + \cdots + a_ps_p = a$
      o\`u $(a_1, \ldots ,a_p)$ est une famille d'\'el\'ements de ${\bf N}$ premier entre eux et $a $ un nombre complexe   est l'\'equation d'un hyperplan affine $H$, les images par 
      ${\rm exp}^{ 2i \pi . }$ des $\tau ^k (H)$ pour $k\in {\bf Z}$ 
co\"{\i}ncident et ont pour \'equation : 
 $$ (\sigma _1)^{a_1} \cdots (\sigma _p)^{a_p} = e^{ 2i \pi a}$$ 
 
  Ce r\'esultat r\'epond \`a une  question de N. Budur \cite{Bu} pos\'ee pour le cas particulier
  $M= {\cal O}_X$.

   \subsubsection{Cas o\`u $m$ est une section d'un Module holonome r\'egulier}

    Nous supposons maintenant que $M = {\cal D}_X m $ est un module holonome r\'egulier de vari\'et\'e caract\'eristique $\Lambda$.
    Notons   $F$ le produit $f_1  \cdots f_p$ et
d\'esignons par $W^{\sharp}_{f_1 , \ldots f_p , \Lambda}$ l'adh\'erence dans $T^{\ast}X \times {\bf C}^p$ de
$$ \Omega = \{ (x, \xi + \sum_{i=1}^p s_i \frac{ df_i(x)}{f_i(x)}, s_1, \ldots ,s_p ) \; ; \; s_i \in {\bf C} \; , \;  (x, \xi) \in  \Lambda 
\; {\rm et} \; F(x) \neq 0 \}\; .$$

Nous avions \'etabli les r\'esultats suivants avec J. Brian\c{c}on et M. Merle    dans \cite{B-M-M1}~:

 \begin{enumerate}
 \item Les fibres r\'eduites de la restriction de $\pi _2$ \`a $ W^{\sharp}_{f_1, \ldots ,f_p, \Lambda}$ sont des sous-espaces lagrangiens de 
 $T^{\ast}X$. La fibre au-dessus de l'origine est   un sous-espace lagrangien conique  not\'ee
 $ W^{\sharp}_{f_1, \ldots ,f_p, \Lambda}(0)$.
 \item Les composantes irr\'eductibles de $ W^{\sharp}_{f_1, \ldots ,f_p, \Lambda} \cap F^{-1}(0)$ sont toutes de dimension $ {\rm dim}\, X +p-1$. Leurs projections par $\pi _2$ sont des  hyperplans vectoriels de ${\bf C}^p$ dont les \'equations sont des formes lin\'eaires \`a coefficients entiers positifs ou nuls. Plus pr\'ecisement, si  $G$ une composante irr\'eductible de $ W^{\sharp}_{f_1, \ldots ,f_p, \Lambda} \cap F^{-1}(0)$, si $n_j$ d\'esigne la multiplicit\'e de $f_j$ le long de $G$, $\pi _2 (G)$
 est l'hyperplan de ${\bf C}^p$ d'\'equation~: $n_1s_1 + \cdots + n_ps_p =0$.
\end{enumerate}
  
  \begin{definition} Nous appelons pentes de $(\Lambda, f_1, \ldots ,f_p )$ les hyperplans vectoriels obtenus par projection par $\pi _2$  des composantes irr\'eductibles de
  $ W^{\sharp}_{f_1, \ldots ,f_p, \Lambda}  \cap F^{-1}(0)  $.  Pou tout $x_0 \in X$, nous notons ${\cal H}(\Lambda,x_0,f_1, \ldots ,f_p )$ l'ensemble de ces pentes  au voisinage de $x_0$.
  \end{definition}

Avec l'aide d'un r\'esultat  de C. Sabbah sur les vari\'et\'e caract\'eristiques  d'un module relatif engendrant un module holonome 
r\'egulier  (th\'eor\`eme 3.2, \cite{S2}),   nous avions   \'etabli avec J. Brian\c{c}on et M. Merle  dans \cite{B-M-M3}  les r\'esultats suivants~:

\begin{proposition} \label{prs} (voir th\'eor\`eme 2.1, 2.4 et 2.6, \cite{B-M-M3}) Soit $M$ un 
 ${\cal D}_X$-Module holonome r\'egulier $M$ de vari\'et\'e caract\'eristique $\Lambda$
  et $m$ une section de $M$ engendrant $M$.  Alors, 
  \begin{enumerate}
\item  $\displaystyle {\rm car}^{\sharp} \, ( {\cal D}_X[s_1, \ldots ,s_p]  m   f_1^{s_1} \ldots  f_p^{s_p} )  =  W^{\sharp}_{f_1, \ldots ,f_p, \Lambda}  \; ,$ 
\item $\displaystyle {\rm car}^{\sharp}\,  \left( \frac{{\cal D}_X[s_1, \ldots ,s_p]  m   f_1^{s_1} \ldots  f_p^{s_p}}{{\cal D}_X[s_1, \ldots ,s_p]  m   f_1^{s_1+1} \ldots  f_p^{s_p+1}} \right)  =  W^{\sharp}_{f_1, \ldots ,f_p, \Lambda}  \bigcap F^{-1}(0)\; .$
\item pour tout $c= (c_1, \ldots,c_p)\in {\bf C}^p$ :
$$ {\rm car}^{\sharp}  \,\left( \frac{{\cal D}_X[s_1, \ldots ,s_p]  m   f_1^{s_1} \ldots  f_p^{s_p}}{ \sum_{j=1}^p (s_j - c_j) {\cal D}_X[s_1, \ldots ,s_p]  m   f_1^{s_1+1} \ldots  f_p^{s_p+1}} \right)  =  W^{\sharp}_{f_1, \ldots ,f_p, \Lambda}(0)\; .$$
\end{enumerate}

\end{proposition}

Nous d\'eduisons des   propositions \ref{pcrf}, \ref{pcvcf} et \ref{pefml}

 \begin{proposition}  Soit $M  = {\cal D}_X m $  un module holonome r\'egulier de vari\'et\'e caract\'eristique $\Lambda$.
 $$  {\rm car}^{\rm rel} \, ( {\cal D}_X[s_1, \ldots ,s_p]  m   f_1^{s_1} \ldots  f_p^{s_p} )  =  W^{\sharp}_{f_1, \ldots ,f_p, \Lambda}(0) \times
 {\bf C}^p \; .$$ 
 \end{proposition}

\begin{proposition} \label{ppbw} Soit $M  = {\cal D}_X m $  un module holonome r\'egulier de vari\'et\'e caract\'eristique $\Lambda$. Il existe un 
polyn\^ome de Bernstein de $m,f_1, \ldots ,f_p$ au voisinage de $x_0$ non nul qui soit produit de formes lin\'eaires
dont l'ensemble des directions vectorielles est exactement   ${\cal H}(\Lambda,x_0,f_1, \ldots ,f_p )$.
\end{proposition}

Cette proposition   avait \'et\'e obtenu pour $p=2$  dans  \cite{B-M-M3}  th\'eor\`eme 3.2.\\
   
   \noindent { \bf Preuve :} Suivant C. Sabbah, consid\'erons un $b_S \in {\cal B}(m,x_0, f_1, \ldots ,f_p)$  non nul produit  de formes lin\'eaires affines \`a coefficients rationnels. Ecrivons $b_S = c_1c_2$ o\`u $c_1$ (resp. $c_2$) est produit de formes lin\'eaires affines de directions vectorielles
   form\'ees par de  pentes de ${\cal H}(\Lambda,x_0,f_1, \ldots ,f_p )$ (resp. non form\'ees).  Nous avons vu qu'au voisinage de $x_0$~:
$b_S(s  -r)\cdots b_S(s  ) \tilde{L}  \subset \tau (   \tilde{L})$.
  Nous en d\'eduisons que le sous-Module~:
  $$  c_1(s  -r)\cdots c_1(s  ) \frac{\tilde{L}}{ \tau (   \tilde{L}) } \subset \tilde{L}  / \tau (   \tilde{L}) $$
 est annul\'e par $  c_2(s -r)\cdots c_2(s  )$. Il en r\'esulte~:
   $$ {\rm car}^{\sharp}  \, (  c_1(s-r)\cdots c_1(s ) \frac{\tilde{L}}{ \tau (   \tilde{L}) })  \subset
   W^{\sharp}_{f_1, \ldots ,f_p, \Lambda}  \bigcap F^{-1}(0)\bigcap ({\rm in } \, c_2)^{-1}(0) $$qui est de
   dimension strictement inf\'erieure \`a ${\rm dim}\, X + p-2$. Comme $\tilde{L}  / \tau (   \tilde{L})$ est pur de grade ${\rm dim}\, X +1 $, il s'en suit~:
   $$ c_1(s_1-r, \ldots ,s_p -r)\cdots c_1(s_1, \ldots ,s_p ) \frac{\tilde{L}}{ \tau (   \tilde{L}) } = 0\; .$$ 
   Posons $b_1(s  )=  c_1(s  -r)\cdots c_1(s  )$, nous obtenons~: 
   $$ b_1(s  +r)\cdots b_1(s  )  \in  {\cal B}(m,x_0, f_1, \ldots ,f_p)\; .$$  
   Comme $ b_1(s  +r)\cdots b_1(s  ) $ est produit de formes lin\'eaires affines de directions vectorielles
   form\'ees par des pentes de ${\cal H}(\Lambda,x_0,f_1, \ldots ,f_p )$, nous avons montr\'e l'existence dans $ {\cal B}(m,x_0, f_1, \ldots ,f_p)$
   d'un tel polyn\^ome non nul. Inversement, la partie homog\`eme de plus haut degr\'e de ce polyn\^ome s'annule sur $W^{\sharp}_{f_1, \ldots ,f_p, \Lambda}  \bigcap F^{-1}(0)$ au voisinage
   de $x_0$ et donc sur toutes les pentes de ${\cal H}(\Lambda,x_0,f_1, \ldots ,f_p )$.

  Si  $c(s_1, \ldots, s_p)$ est un polyn\^ome, d\'esignons par ${\rm in }\, c$ la partie homog\`ene de $c$ de plus haut degr\'e. Si $b$ appartient \`a l'ideal de Bernstein de $m,f_1, \ldots ,f_p$ au voisinage de $x_0$, ${\rm in }\, b$ s'annule sur 
  $$ {\rm car}^{\sharp} \, \left( \frac{{\cal D}_X[s_1, \ldots ,s_p]  m   f_1^{s_1} \ldots  f_p^{s_p}}{{\cal D}_X[s_1, \ldots ,s_p]  m   f_1^{s_1+1} \ldots  f_p^{s_p+1}} \right)\; . $$
Il r\'esulte de la proposition \ref{prs} que ${\rm in }\, b$ s'annule sur  ${\cal H}(\Lambda,x_0,f_1, \ldots ,f_p )$.\\

La proposition \ref{ppbw}  permet de montrer que  les  pentes ${\cal H}(x_0,m)$  de $(m,f_1, \dots,f_p)$ au voisinage de $x_0$ 
ne sont autres que 
  ${\cal H}(\Lambda,x_0,f_1, \ldots ,f_p )$ l'ensemble de ces pentes  au voisinage de $x_0$. Ainsi, nous obtenons :


\begin{proposition} \label{psibr}Soit $M  = {\cal D}_X m $  un module holonome r\'egulier de vari\'et\'e caract\'eristique $\Lambda$.
Au voisinage de tout $x_0  \in  X$, il existe 
      $(X_{\alpha})_{\alpha \in A}$ (resp. $S_{\alpha}$) une  famille finie de sous-espaces analytiques (resp. vari\'et\'es alg\'ebriques de ${\bf C}^p$)
        
      telle  que 
$$  {\rm car}^{\rm rel }\,  \left(  \frac{{\cal D}_X[s_1, \ldots ,s_p] m f_1^{s_1} \ldots  f_p^{s_p}}{{\cal D}_X[s_1, \ldots ,s_p] m f_1^{s_1+1} \ldots  f_p^{s_p+1}}  \right)  =  \cup_{\alpha \in A}T^{\ast}_{X_{\alpha}}X \times  S_{\alpha}   \; .$$
avec les conditions
 \begin{itemize}
 \item Chaque vari\'et\'e alg\'ebrique  $S_{\alpha}$ est de dimension $p-1$.  
 \item Les composantes irr\'eductibles de 
 dimension $p-1$ de chaque $S_{\alpha}$  sont des hyperplans affines  $H_{\alpha,\beta}$  de directions les noyaux de  formes lin\'eaires \`a coefficients premiers entre eux dans ${\bf N}$.
 \item   Toute   composante  irr\'eductible   des  $S_{\alpha}$  de    dimension strictement inf\'erieure \`a $p-1$ est  contenue  dans un hyperplan affine $\tau ^k ( H_{\alpha,\beta}  )$ o\`u $k \in {\bf Z}$ et 
  $\tau$ la translation $(s_1, \ldots, s_p) \mapsto (s_1+1, \ldots, s_p+1)$. 
  \item  l'ensemble des  directions des hyperplans affines $H_{\alpha,\beta}$ pour $x_0 \in X_{\alpha}$ est l'ensemble  ${\cal H}(\Lambda,x_0,f_1, \ldots ,f_p )$.\\
  \end{itemize} 
 \end{proposition}

 \begin{corollaire} Soit $M  = {\cal D}_X m $  un module holonome r\'egulier de vari\'et\'e caract\'eristique $\Lambda$. La vari\'et\'e des z\'eros de    
  ${\bf in }\,  {\cal B}(m,x_0, f_1, \ldots ,f_p)$  l'id\'eal   engendr\'e par les  parties homog\'enes de plus haut degr\'e  
     des \'el\'ements de $ {\cal B}(m,x_0, f_1, \ldots ,f_p)$ est la reunion est la r\'eunion des pentes de $(\Lambda, f_1, \ldots ,f_p )$ au voisinage de $x_0$.
  
         \end{corollaire}
         
          \noindent { \bf Preuve :} Si $b\in  {\cal B}(m,x_0, f_1, \ldots ,f_p)$, $ {\rm in}\,  b$ la partie homog\`ene de plus haut  degr\'e de $b$ s'annule sur la vari\'et\'e caract\'eristique di\`ese de 
          $$ {\cal D}_X[s_1, \ldots ,s_p]  m   f_1^{s_1} \ldots  f_p^{s_p}/{\cal D}_X[s_1, \ldots ,s_p]  m   f_1^{s_1+1} \ldots  f_p^{s_p+1} $$
          et donc, suivant la proposition \ref{prs}, s'annule sur 
 $ W^{\sharp}_{f_1, \ldots ,f_p, \Lambda}  \bigcap F^{-1}(0) $. Il en r\'esulte une inclusion. L'autre inclusion est donn\'ee par la proposition \ref{ppbw}.\\

 Enfin, nous avons :

 \begin{corollaire} Soit $M  = {\cal D}_X m $  un module holonome r\'egulier de vari\'et\'e caract\'eristique $\Lambda$. L'image par l'application ${\rm exp}^{ 2i \pi . }$ de la  vari\'et\'es des z\'eros de $  {\cal B}(m,x_0, f_1, \ldots ,f_p)$ est une r\'eunion
        de sous-ensembles de $ ({\bf C}^{\ast})^p$  o\`u  chaque sous-ensemble est d\'efini  par une \'equation dun type
  $$ (\sigma _1)^{a_1} \cdots (\sigma _p)^{a_p} = \alpha $$ 
  o\` u   $\alpha $ est un nombre complexe et    $(a_1, \ldots ,a_p)$ $(a_1, \ldots ,a_p)$ est une famille d'\'el\'ements de ${\bf N}$ premier entre eux. Cette famille est 
exactement  la famille  des coefficients des formes lin\'eaires dont les des z\'eros sont 
  les pentes de $(\Lambda, f_1, \ldots ,f_p )$ au voisinage de $x_0$.
  
 \end{corollaire}

\end{document}